\documentclass[10pt,twoside]{article}
\usepackage{amssymb}
\usepackage{amsmath}

\newtheorem{thm}{Theorem}[section]

\newtheorem{propn}[thm]{Proposition}
       
\newtheorem{defn}[thm]{Definition}

\newtheorem{lmm}[thm]{Lemma}

\newtheorem{cor}[thm]{Corollary}

\newtheorem{remark}[thm]{Remark}

\newcommand{\R}{{\bf R}}
\newcommand{\Z}{{\bf Z}}

\newcommand{\C}{{\bf C}}
\newcommand{\Q}{{\bf Q}}
\newcommand{\ad}{\overline{\mbox{ad}}}

\newenvironment{Pf}{\medskip \noindent {\bf Proof: }}
   {\hfill $\diamondsuit$ }

\begin{document}

\title{Compact Lorentz manifolds with local symmetry}
\author{Karin Melnick}
\date{\today}
\maketitle

\newpage

\section{Acknowledgements}

First I would like to thank my advisor, Benson Farb, for his guidance, inspiration, and encouragement.  Working with him has been a great pleasure for which I am deeply grateful.  He is a remarkable source of ideas, and I benefitted from his knowledge of a broad range of mathematics.

While working on this project, I had many helpful conversations with Abdelghani Zeghib.  He suggested approaching this problem with Gromov's theory of rigid geometric structures.  His results on codimension-one, totally geodesic, lightlike foliations, among others, play an important role in this work.  

I was very fortunate to have been exposed to the work of Benson Farb and Shmuel Weinberger on compact aspherical Riemannian manifolds with local symmetry.  Their results are of course the fundamental inspiration for this project.  I am priveleged to have learned their remarkable new techniques from them.  I am grateful in particular for several helpful conversations with my second advisor, Shmuel Weinberger.  I would also like to thank my topic advisor, Robert Zimmer, who first brought my attention to this work.  All I learned from him during 2001-2002 continues to influence my thinking about mathematics.

Finally, I enjoyed helpful conversations with many others, including Thierry Barbot, Mark Behrens, Bill Goldman, and Ben Wieland.  I am particularly grateful to Thierry Barbot for his comments and corrections on a preliminary draft of this thesis.

\bigskip

\emph{In loving memory of Laura Sue-Jung Kang}

\newpage

\section{Introduction}

Two main lines of questioning in the study of automorphism groups of pseudo-Riemmanian manifolds are to ask which groups can act isometrically on pseudo-Riemannian manifolds of a given topological type, and to ask which pseudo-Riemannian manifolds admit an isometric action of a given large group.  This work addresses a modified version of the second question: which compact Lorentz manifolds admit a positive-dimensional pseudogroup of local isometries?  This question can be loosely rephrased as, which compact Lorentz manifolds have nontrivial local symmetry?  For a real-analytic, compact, complete Lorentz manifold, a positive-dimensional pseudogroup of local isometries is equivalent to a positive-dimensional isometry group on the universal cover.

Several examples of compact Lorentz manifolds with local symmetry will be discussed below.  Given such a Lorentz manifold, one may construct a new compact Lorentz manifold with at least as much local symmetry by forming a warped product.

\begin{defn}
For two pseudo-Riemannian manifolds $(P, \lambda)$ and $(Q,\mu)$, a \emph{warped product} $P \times_f Q$ is given by a positive function $f$ on $Q$ : the metric at $(p,q)$ is $f(q) \lambda_p + \mu_q$.  The factor $P$ is called the \emph{normal factor}.
\end{defn}

If $Isom(P) = G$, then $G$ also acts isometrically on the warped product $P \times_f Q$ for any $f$.  More generally, let $f$ be any function $Q \rightarrow \mathcal{M}$, where $\mathcal{M}$ is the moduli space of $G$-invariant metrics (of a fixed signature) on $P$, with $f(q) = \lambda(q)$.  Then $G$ acts isometrically on the \emph{generalized warped product} $P \times_f Q$, where the metric at $(p,q)$ is $\lambda(q)_p + \mu_q$.

Results of Farb and Weinberger stated below give conditions under which a compact Riemannian manifold is a generalized warped product $P \times_f Q$ with $P$ a locally symmetric space.
Our main result (Theorem \ref{mainthm} below) gives conditions under which the universal cover of a compact Lorentz manifold has this form with $P$ a Riemannian symmetric space or a complete Lorentz space of constant curvature.  In both cases, the conditions are that the manifold have a large pseudogroup of local isometries.

Pseudo-Riemannian metrics are examples of \emph{rigid geometric structures of algebraic type}.  For $M$ a compact real-analytic manifold with such a structure, Gromov's Stratification Theorem (stated as Theorem \ref{stratification} below) describes the orbit structure of local symmetries of $M$.  The celebrated Open-Dense Theorem, which is a corollary of this stratification, states that if a point of $M$ has a dense orbit under local isometries, then an open dense subset of $M$ is locally homogeneous.  It would be interesting to find conditions on $M$ under which existence of a dense orbit implies that $M$ is locally homogeneous. More generally, one might seek a fibered version: when does existence of a local isometry orbit with positive-dimensional closure imply that $M$ is roughly a fiber bundle with locally homogeneous fibers?  Our main theorem (\ref{mainthm}) can be viewed as such a result, under some particular topological and geometric conditions on a compact real-analytic Lorentz manifold.

\subsection{Riemannian case}

For $M$ a compact Riemannian manifold, $Isom(M)$ is compact.  For example, a compact locally symmetric space of noncompact type has finite isometry group (see \cite{wittearith} 5.43).  While such a group provides some information about $M$, the isometry group of the universal cover $X$ of $M$ tells much more.  For example, if $M$ is a locally symmetric space of noncompact type, then $Isom(X)$ is a semisimple group with no compact factors.  A homogeneous, contractible, Riemannian manifold with this isometry group must be a symmetric space.  

Recall that an \emph{aspherical} manifold is one with contractible universal cover.
Farb and Weinberger studied compact aspherical Riemannian manifolds $M$ with universal cover $X$ having $Isom^0(X) \neq 1$.  They proved several results characterizing warped products with locally symmetric factors, and locally symmetric spaces in particular.  The following theorem is a weakened statement of their main theorem.  Orbibundles will be defined later below.

\begin{thm}[Farb and Weinberger \cite{FW}]  Let $M$ be a compact aspherical Riemannian manifold with universal cover $X$. Let $G = Isom(X)$.  If $G^0 \neq 1$, then $M$ is a Riemannian orbibundle
$$ \Lambda \backslash G^0 /K \rightarrow M \rightarrow Q$$

where $\Lambda \subset G^0$ is a cocompact lattice, $K$ is a maximal compact subgroup of $G^0$, and $Q$ is aspherical.  

Further, if $\pi_1(M)$ contains no normal free abelian subgroup, then 
$Z(G^0)$ is finite,  $G^0$ is semisimple, and a finite cover of $M$ is isometric to 
$$\Lambda \backslash G^0 / K \times_f Q$$  

for $f : Q \rightarrow \mathcal{M}$, the moduli space of locally symmetric metrics on $\Lambda \backslash G^0/K$.
\end{thm}

The aspherical assumption is required.  There are metrics on the sphere $S^n$ with a bump at one point, for example, for which the isometry group contains only rotations fixing that point.  However, \cite{FW} contains the statement of a similar theorem to the above, under a noncompactness assumption on the connected isometry group, for arbitrary closed Riemannian manifolds.

Their proof relies on the theory of proper transformation groups, Lie theory, and remarkable cohomological dimension arguments.

\subsection{Lorentz case}

For Lorentz manifolds, a crucial difference from the Riemannian case is that the isometry group need not act properly; in particular, orbits may not be closed, and the group of deck transformations may not intersect $G^0$ in a lattice.  On the other hand, fantastic work has been done on nonproper isometric Lorentz actions (\cite{Kowthesis}, \cite{Zetgl1}, \cite{Zetgl2}), which implies a great deal of structure in that case.  


For a compact Lorentz manifold $M$, the groups $Isom^0(M)$ have been classified (\cite{ZiLor}, \cite{AS1}, \cite{AS2}, \cite{Ze1}, \cite{Ze2}).  There are several results on the form of a compact Lorentz manifold admitting an isometric action of a given group (\cite{Gromov} 5.4.A, \cite{Ze2} 1.14, \cite{Ze2} 4.1.2, \cite{Me} 4.9).  Here we consider universal covers of compact Lorentz manifolds with isometric actions of semisimple groups; some techniques on compact manifolds can be extended to this setting.

The Lorentz manifolds with the most symmetry are those of constant curvature, modelled on Minkowski space, de Sitter space, or anti-de Sitter space.  Any \emph{Lorentzian locally symmetric space} has constant curvature, as was proved in \cite{CLPTV} and independently in \cite{ZeCC}.  A \emph{Lorentzian locally symmetric space} is a Lorentz manifold such that, for any tangent vector ${\bf v}$ at any point $x$, there is a local isometry fixing $x$ and sending ${\bf v}$ to $- {\bf v}$.  The model for \emph{$n$-dimensional Minkowski space}, $Min^n$, is $R^n$ equipped with a nondegenerate inner-product of type $(1,n-1)$.   The model for \emph{$n$-dimensional de Sitter space}, $dS^n$, is the $+1$-level set of a quadratic form of type $(1,n)$, with the induced metric.  The model for \emph{$n$-dimensional anti-de Sitter space}, $AdS^n$, is the $-1$-level set of a quadratic form of type $(2,n-1)$, with the induced metric.  Each is a homogeneous space, $G/H$, where $H$ is the stabilizer of a point.  The isometry group, stabilizer, curvature, and diffeomorphism type for each are in the following table.

\medskip
\begin{center}
\begin{tabular}{|l|c|c|c|}
\hline
       &  $Min^n$                & $dS^n$               & $AdS^n$ \\
\hline
       &                         &                      &          \\
Isom   &  $O(1,n-1) \ltimes R^n$ & $O(1,n)$             & $O(2,n-1)$    \\
       &                         &                      &           \\
Stab   &  $O(1,n-1)$             & $O(1,n-1)$           & $O(1,n-1)$     \\
       &                         &                      &              \\
Curv   &  $0$                    &  $1$                 & $-1$  \\
       &                         &                      &        \\
Diff   &   $\R^n$                & $S^{n-1} \times \R$  & $\R^{n-1} \times S^1$ \\
\hline
\end{tabular}
\end{center}
\medskip

Note that $AdS^2 \cong dS^2 \cong SO(1,2) / A$, where $A \cong \R^*$ is a maximal $\R$-split torus.  On Minkowski space of any dimension, there are obviously discrete groups of isometries acting properly discontinuously and cocompactly.  A result of Calabi and Markus states that no infinite subgroup of $O(1,n)$ acts properly on $dS^n$, so there are no compact complete de Sitter manifolds (\cite{CM}).  Kulkarni noted that when $n$ is odd, lattices in $SU(1,(n-1)/2)$ act freely, properly discontinously, and cocompactly on $AdS^n$.  For $n$ even, he proved that there is no cocompact, properly discontinuous, isometric action on $AdS^n$ (\cite{Kulk}).  The group $SL_2(\R)$ with the Killing metric is isometric to $AdS^3$.   

Kowalsky, using powerful dynamical techniques, which are treated in detail in Section \ref{kowalsky.argument} below, proved that a simple group acting nonproperly on an arbitrary Lorentz manifold is locally isomorphic to $O(1,n)$, $n \geq 2$, or $O(2,n)$, $n \geq 3$ (\cite{Kowthesis}).  Adams has characterized groups that admit \emph{orbit nonproper} isometric actions on arbitrary Lorentz manifolds in \cite{Adamsonp1} and \cite{Adamsonp2}; an action $G \times M \rightarrow M$ is \emph{orbit nonproper} if for some $x \in M$, the map $g \mapsto g.x$ from $G$ to $M$ is not proper.

There are several recent results on the form of arbitrary Lorentz manifolds admitting isometric actions of certain semisimple groups.  Witte Morris showed that a homogeneous Lorentz manifold with isometry group $O(1,n)$ or $O(2,n-1)$ is $dS^n$ or $AdS^n$, respectively (\cite{WM}).  Arouche, Deffaf, and Zeghib, using totally geodesic, lightlike hypersurfaces, showed that if a semisimple group with no local $SL_2(\R)$-factors has a Lorentz orbit with noncompact isotropy, then a neighborhood of this orbit is a warped product $N \times_f L$, where $N$ is a complete, constant-curvature Lorentz space, and $L$ is a Riemannian manifold (\cite{ADZ}).  Deffaf, Zeghib, and the author treat degenerate orbits with noncompact isotropy in \cite{DMZ}.  We conclude that any nonproper action of a semisimple group with finite center and no local $SL_2(\R)$-factors has an open subset isometric to a warped product as in \cite{ADZ}, and we describe the global structure of such actions.

The work here combines features and techniques of many of these papers, as well as those of \cite{FW}.  As in \cite{FW}, we consider universal covers of compact aspherical Lorentz manifolds and seek to describe those for which the identity component of the isometry group is nontrivial.  Here is the main result.

\begin{thm}  
\label{mainthm}
Let $M$ be a compact, aspherical, real-analytic, complete Lorentz manifold with universal cover $X$.  Let $G = Isom(X)$, and assume $G^0$ is semisimple.  

{\bf (1)  Orbibundle.}  Then $M$ is an orbibundle

$$ P \rightarrow M \rightarrow Q$$
where $P$ is aspherical and locally homogeneous, and $Q$ is a good aspherical orbifold.

{\bf (2)  Splitting.} Further, precisely one of the following holds:

{\bf A. } $G^0$ acts properly on $X$: 

Then $P = \Lambda \backslash G^0 / K$ where $\Lambda$ is a lattice in $G^0$ and $K$ is a maximal compact subgroup of $G^0$.

Further, if $|Z(G^0)| < \infty$, then a finite cover of $M$ is isometric to 

$$P \times_f Q$$  

for $f : Q \rightarrow \mathcal{M}$, the moduli space of Riemannian locally symmetric metrics on $P = \Lambda \backslash G^0/K$.  The Lorentzian manifold $Q$ has $Isom^0(\widetilde{Q}) = 1$.

{\bf B. } $G^0$ acts nonproperly on $X$: 

Then $M$ is a Lorentzian orbibundle.  The metric along $G^0$-orbits is Lorentzian, with
$$P = \Lambda \backslash (\widetilde{AdS}^k \times G_2 / K_2)$$ 

where $k \geq 3$, $G_2 \lhd G^0$ with maximal compact subgroup $K_2$, and 
$$\Lambda \subset \widetilde{O}^0(2,k-1) \times G_2$$ 

acts freely, properly discontinuously, and cocompactly on $\widetilde{AdS}^k \times G_2/K_2$.  The good Riemannian orbifold $Q$ has $Isom^0(\widetilde{Q}) = 1$.  There is a warped product

$$X \cong  \widetilde{AdS}^{k} \times_h L$$

for some real-analytic function $h : L \rightarrow \R^+$.

Further, if $|Z(G_2)| < \infty$, then $X$ is isometric to 
$$ (\widetilde{AdS}^k \times G_2 / K_2) \times_f \widetilde{Q}$$

where $f : \widetilde{Q} \rightarrow \mathcal{M}$, and $\mathcal{M} \cong \R^2$ is the moduli space of $G^0$-invariant Lorentzian metrics on $\widetilde{AdS}^k \times G_2/K_2$.
\end{thm}

\begin{cor}
Let $M$ and $G^0$ be as above.  If $M$ has an open, dense, locally homogeneous subset, then $M$ is locally homogeneous.
\end{cor}

Section \ref{section.examples} below contains a construction in which a noncompact, connected, semisimple group $H^0 \subseteq Isom(X)$; $Z(H^0)$ is infinite; $H^0$ acts properly on $X$; and the metric type on $H^0$-orbits varies, which we believe illustrates the necessity of the hypothesis of finite center in (2) {\bf A}.  In that section, we also adapt to the Lorentz setting a construction of \cite{FW} of \emph{essential orbibundles}---that is, orbibundles that are not finitely covered by any fiber bundle.

{\bf Proof Outline} for Theorem \ref{mainthm}:

\begin{itemize}

\item{The first step involves Gromov's stratification for isometric actions on spaces with rigid geometric structure: there is a closed orbit in $X$ on which the group of deck transformations acts cocompactly (Propositions \ref{goodorbit}, \ref{gamma0.cocompact}). }

\end{itemize}

The stabilizer of a point in this orbit then determines the dynamics of the isometry group on $X$.

\begin{itemize}

\item{If the stabilizer is compact, then the group generated by $G^0$ and the fundamental group acts properly.  In this case, techniques of \cite{FW} apply (Section \ref{section.proper.case}).}

\item{When the stabilizer is noncompact, then $G^0$ acts nonproperly on $X$.  In this case, we extend work of \cite{Zetgl1} to show that totally geodesic lightlike foliations exist on $X$ (Theorem \ref{tglsexist}).  When there are many of these foliations, then results of \cite{Zetgl2} give the warped product structure on $X$ and the orbibundle in $M$ (Sections \ref{section.warped.prod}, \ref{section.orbibundle}).  We show in Section \ref{tglx.big} that there must be sufficiently many of these foliations; the argument here involves dynamical techniques of \cite{Kowthesis}, Lie algebras, and some basic facts about locally homogeneous spaces.}
\end{itemize}

\section{Notation}

Throughout, $M$ is a compact, aspherical, real-analytic, complete Lorentz manifold.  The universal cover of $M$ is $X$, with $Isom(X) = G.$   The group of deck transformations is $\Gamma \cong \pi_1(M)$.  The identity component of $G$ is a semisimple group $G^0$, and $\Gamma_0 = \Gamma \cap G^0$.  Note $G^0 \lhd G$ and $\Gamma_0 \lhd \Gamma$.

The Lie algebra of $G^0$ is $\mathfrak{g}$.  Let $\mathfrak{g} = \mathfrak{g}_1 \oplus \cdots \oplus \mathfrak{g}_l$ be the decomposition of $\mathfrak{g}$ into simple factors.  Let $G_i$ be the corresponding subgroups of $G^0$.  The projection $\mathfrak{g} \rightarrow \mathfrak{g}_i$ will be denoted $\pi_i$, as will the projection $G^0 \rightarrow G_i$.  

For an arbitrary group $H$ acting on a space $Y$, the stabilizer of $y \in Y$ will be denoted $H(y)$.  In particular, $G_i(y) = G^0(y) \cap G_i$, and $\mathfrak{g}_i(y) = \mathfrak{g}(y) \cap \mathfrak{g}_i$.

\section{Background and Terminology}

\subsection{Proper actions}

Let $H$ be a locally compact topological group and $Y$ a locally compact Hausdorff space.

\begin{defn}  A continuous action of $H$ on $Y$ is \emph{proper} if, for any compact subsets $A,B \subset Y$, the set 
$$ H_{A,B} = \{ h : hA \cap B \neq \emptyset \}$$

is compact in $H$.  
\end{defn}

Note that $H_{A,B}$ is automatically closed by continuity of the action.  The following equivalence is easy to show.  Both characterizations of properness will be used below.

\begin{propn}
The action of $H$ on $Y$ as above is proper if and only if, for any compact $A \subseteq Y$, the set
$$ H_A = \{ h : hA \cap A \neq \emptyset \}$$

is compact in $H$.
\end{propn}

\begin{defn} A proper action of a discrete group is called \emph{properly discontinuous}.
\end{defn}

The following facts about proper actions are standard.

\begin{propn}
If $H$ acts properly on $Y$, then 
\begin{enumerate}
\item{$H(y)$ is compact for all $y \in Y$}
\item{$Hy$ is closed in $Y$}
\item{$H \backslash Y$ is Hausdorff in the quotient topology}
\item{Any closed subgroup of $H$ also acts properly on $Y$}
\end{enumerate}
\end{propn}

The following facts rely on the existence of \emph{slices} for smooth proper actions on manifolds. See \cite{Pflaum} for definitions related to stratified spaces.

\begin{propn}
\label{stratfacts}
Let $H$ be a Lie group acting smoothly and properly on a connected manifold $Y$.
\begin{enumerate}
\item{ For any compact $\overline{A} \subset H \backslash Y$, there is a compact $A \subset Y$ projecting onto $\overline{A}$.}

\item{ In general, $H \backslash Y$ is a Whitney stratified space with 
$$\mbox{dim} (H \backslash Y) = \mbox{dim} Y - \mbox{dim} H + \mbox{dim} H(y)$$

where $\mbox{dim} H(y)$ is minimal over $y \in Y$.}

\item{If the stabilizers $H(y)$ belong to the same conjugacy class for all $y \in Y$, then $H \backslash Y$ is a smooth manifold.}

\end{enumerate}
\end{propn}

\begin{Pf}
\begin{enumerate}
\item{Let $\pi$ be the projection $Y \rightarrow H \backslash Y$.  For any $\overline{y} = \pi(y)$ in $ H \backslash Y$, the Slice Theorem (see \cite{Pflaum} 4.2.6) gives a neighborhood $\overline{U}$ of $\overline{y}$ and a diffeomorphism $\varphi_y :  H \times_{H(y)} V_y \rightarrow \pi^{-1}(\overline{U})$, where $V_y$ is an open ball in some $\R^k$. A disk about ${\bf 0}$ in $V_y$ corresponds under $\varphi_y$ to a compact $D_y$ containing $y$, and projecting to a compact neighborhood of $\overline{y}$ in $H \backslash Y$.  For a compact subset $\overline{A}$, there exist $\overline{y}_1, \ldots, \overline{y}_n$ such that $int(\pi(D_{y_1})), \ldots, int(\pi(D_{y_n}))$ cover $\overline{A}$.  Then $D_{y_1} \cup \cdots \cup D_{y_n}$ is the desired compact $A \subset Y$.  }

\item{The stratification is by orbit types: for each compact $K \subset H$, let 
$$Y_{(K)} = \{ y \in Y \ : \ g H(y) g^{-1} = K \ \mbox{for some}\ g \in H \}$$

 and let $Y_K$ be the fixed set of $K$.  Then the pieces of the stratification of $ H \backslash Y$ are the components of the quotients 
$$ H \backslash Y_{(K)}  = N_H (K) \backslash Y_K$$

Each piece has the structure of a smooth manifold.  See \cite{Pflaum} 4.3.11 and 4.4.6.  When $K = H(y)$ is minimal, then $Y_{(K)}$ is open, and the piece $H \backslash Y_{(K)}$ has maximal dimension $\mbox{dim} Y - \mbox{dim}H + \mbox{dim} H(y)$.}

\item{If all stabilizers are conjugate to one compact subgroup $K$, then $ H \backslash Y = H \backslash Y_{(K)}$, which consists of a single piece, because $Y$ is connected.}
\end{enumerate}
\end{Pf}

\begin{thm}[Goresky \cite{Goresky}]
\label{strattriang}
Stratified spaces can be triangulated such that the interior of each simplex is contained in a piece of the stratification.
\end{thm}


\subsection{Orbifolds and orbibundles}

\begin{defn} An \emph{$n$-dimensional orbifold} is a Hausdorff, paracompact space with an open cover $\{ U_i \}$, closed under finite intersections, with homeomorphisms
$$\varphi_i : \widetilde{U}_i / \Lambda_i \rightarrow U_i$$

where $\widetilde{U}_i$ is an open subset of $\R^n$ and $\Lambda_i$ is a finite group.  The atlas $(U_i, \varphi_i)$ must additionally satisfy the compatibility condition: whenever $U_j \subset U_i$, then there is a monomorphism $\Lambda_j \rightarrow \Lambda_i$ and an equivariant embedding $\widetilde{U}_j \rightarrow \widetilde{U}_i$ inducing a commutative diagram
$$
\begin{array}{ccc}
\widetilde{U}_j & \rightarrow & \widetilde{U}_i \\
\downarrow &                  & \downarrow \\
\widetilde{U}_j/\Lambda_j & \rightarrow & \widetilde{U}_i / \Lambda_i \\
\downarrow &                  & \downarrow \\
U_j & \rightarrow & U_i
\end{array}
$$

A \emph{smooth orbifold} is an orbifold for which the action of each $\Lambda_i$ is smooth, and the embeddings $\widetilde{U}_j \rightarrow \widetilde{U}_i$ are smooth.
\end{defn}

\begin{defn}  A \emph{good (pseudo-Riemannian) orbifold} is the quotient of a (pseudo-Riemannian) manifold by a smooth, properly discontinuous, (isometric) action.  
\end{defn}

It is not hard to see using proper discontinuity that a good orbifold is a smooth orbifold. 

\begin{defn} A smooth \emph{orbibundle} is a manifold $M$ with a projection $\pi$ to a good orbifold $B$, written
$$N \rightarrow M \stackrel{\pi}{\rightarrow} B$$

where $N$ is a manifold, and the orbifold charts $(U_i,\varphi_i)$ on $B$ lift to
$$\psi_i : \pi^{-1}(U_i) \rightarrow N \times_{\Lambda_i} \widetilde{U_i}$$

where $\Lambda_i$ acts freely and smoothly on $N \times \widetilde{U}_i$.  

A \emph{pseudo-Riemannian orbibundle} is a pseudo-Riemannian manifold $M$ with a projection $\pi$ to a good pseudo-Riemannian orbifold $B$ as above, such that the maps arising from $\psi_i$ 
$$\widetilde{U}_i \rightarrow N \times_{\Lambda_i} \widetilde{U}_i \rightarrow \pi^{-1}(U_i)$$

 are isometric immersions.
\end{defn}

Note that, for a pseudo-Riemannian orbibundle $M$ the type of the metric on $M$ may be different from the type of the metric on the quotient orbifold $B$.

\subsection{Rational cohomological dimension}

\begin{defn}  The \emph{integral cohomological dimension} of a group $\Lambda$ is
$$\mbox{cd}_{\Z} \Lambda = sup \{n\ :\ H^n(\Lambda, A) \neq 0, \ A \ \mbox{a $\Lambda \Z$-module} \}$$

The \emph{rational cohomological dimension} of $\Lambda$ is
$$\mbox{cd}_{\Q} \Lambda = sup \{n\ :\ H^n(\Lambda,A) \neq 0, \ A \ \mbox{a $\Lambda \Q$-module} \}$$
\end{defn}

Proofs of the following facts about integral cohomological dimension can be found in \cite{Brown}; exact references are given for each item.

\begin{propn}  
\label{basiccdz}
Let $\Lambda$ be a discrete group. 
\begin{enumerate}
\item{(VIII.2.2) Let $Y$ be a contractible space on which $\Lambda$ acts freely and properly with $\Lambda \backslash Y$ a finite CW-complex.  Then 
$$\mbox{cd}_{\Z} \Lambda \leq \mbox{dim} Y$$
}

\item{(VIII.2.4b)  Let $\Lambda_0 \lhd \Lambda$.  Then
$$\mbox{cd}_{\Z} \Lambda \leq \mbox{cd}_{\Z} \Lambda_0 + \mbox{cd}_{\Z} (\Lambda / \Lambda_0)$$
.}

\item{(VIII.2.5) If $\mbox{cd}_{\Z} \Lambda$ is finite, then $\Lambda$ is torsion-free.}
\end{enumerate}
\end{propn}

Next we collect some facts about rational cohomological dimension.

\begin{propn}
\label{basiccdq}
Let $\Lambda$ be a discrete group.
\begin{enumerate}

\item{Let $Y$ be a contractible space on which $\Lambda$ acts freely and properly with $\Lambda \backslash Y$ a finite CW-complex.  Then
$$\mbox{cd}_{\Q} \Lambda \leq \mbox{dim} Y$$

If $\Lambda \backslash Y$ is a manifold, then there is equality.
 }

\item{If $\Lambda$ is finite, then $\mbox{cd}_{\Q} \Lambda = 0$.}

\item{Let $\Lambda_0 \lhd \Lambda$.  Then
$$\mbox{cd}_{\Q} \Lambda \leq \mbox{cd}_{\Q} \Lambda_0 + \mbox{cd}_{\Q} (\Lambda / \Lambda_0)$$
}

\end{enumerate}
\end{propn}

\begin{Pf}
\begin{enumerate}

\item{Under these hypotheses, $\Lambda \backslash Y$ is a $K(\Lambda, 1)$ space.  Then $H^*(\Lambda, A) \cong H^*(\Lambda \backslash Y, A)$ for any $\Lambda$-module $A$.  For any $A$, the cohomology $H^n( \Lambda \backslash Y , A) = 0$ for $n > \mbox{dim} Y$, yielding the desired inequality.  If the quotient space is an orientable manifold, then $H^n(\Lambda \backslash Y, \Z) \neq 0$ when $n = \mbox{dim} Y$.
}

\item{This statement follows from (1).}

\item{This inequality is derived from the Hochschild-Serre spectral sequence, (see \cite{Brown} VII.6.3), which applies to group cohomology with coefficients in arbitrary modules.}
\end{enumerate}
\end{Pf}

Rational cohomological dimension gives information about actions that are properly discontinuous but not necessarily free.  The following fact about such actions will be cited multiple times below.

\begin{propn}
\label{cellular}
Let $\Lambda$ be a discrete group acting properly discontinuously on a locally compact Hausdorff topological space $Y$ such that the quotient $\Lambda \backslash Y$ admits a CW decomposition compatible with the stratification by orbit types---the interior of each cell is contained in a piece of the stratification. Then $Y$ admits a CW decomposition on which the $\Lambda$-action is cellular.
\end{propn}

\begin{Pf}
(compare \cite{Schubert} III.6.9.2)
Denote by $\pi$ the projection $Y \rightarrow \Lambda \backslash Y$.
For a finite subgroup $F$ of $\Lambda$, denote by $Y_{(F)}$ the stratum consisting of all points $y$ for which $\Lambda(y)$ is conjugate to $F$.  The set $Y_{(F)}$ is $\Lambda$-invariant.  Let $\overline{Y}_{(F)} = \pi(Y_{(F)})$.  The collection of $\overline{Y}_{(F)}$ form the stratification by orbit types of $\Lambda \backslash Y$. 

We will show that the restriction of $\pi$ to any stratum in $Y$ is a covering map.  Let $x \in \overline{Y}_{(F)}$ and $y \in \pi^{-1}(x)$.  To find a neighborhood $U$ of $y$ in $Y_{(F)}$ that projects homeomorphically to its image under $\pi$, it suffices to find a neighborhood on which $\pi$ is injective, because $\pi$ is continuous and open.  Assume that $\Lambda(y) = F$.  By properness, there exists a neighborhood $U$ of $y$ in $Y$ such that $\lambda U  \cap U = \emptyset$ for all $\lambda \notin F$.  For any $z \in U$, the stabilizer $\Lambda(z) \subseteq F$.  If $\pi(z) = \pi(w)$ for some $z,w \in U$, then $z = \lambda w$ for some $\lambda \in F$.  Now for every $z \in U \cap Y_{(F)}$, we have $\Lambda(z) = F$; therefore, $\pi$ is injective on this neighborhood in $Y_{(F)}$.  Let $V = \cap_{\lambda \in F} \lambda( U \cap Y_{(F)})$.  Then $\pi^{-1}(\pi(V)) = \Lambda V$ is a disjoint collection of neighborhoods in $Y_{(F)}$, each mapping homeomorphically to its image in $\overline{Y}_{(F)}$.  

For each $0$-cell $\overline{e}$ of $\Lambda \backslash Y$, choose a lift $e_{1}$ to $Y$, and let $e_{\lambda} = \lambda e_{1}$ be $0$-cells of $Y$, where $\lambda$ ranges over a collection of coset representatives for $\Lambda / \Lambda(e_1)$.  Now suppose that the $(k -1)$-skeleton of $Y$ has been constructed, and that it is the inverse image of the $(k-1)$-skeleton in $\Lambda \backslash Y$.  Denote by $D^k$ the $k$-dimensional disk, and by $B^k$ its interior.  Given a $k$-cell $\overline{e}$ with $f_{\overline{e}} : D^k \rightarrow \Lambda \backslash Y$, say $f_{\overline{e}}(B^k) \subseteq \overline{Y}_{(F)}$.  Because $B^k$ is contractible, there is a lift $f_{e_1} : B^k \rightarrow Y_{(F)}$.  This lift is a homeomorphism onto its image.  The extension of $f_{e_1}$ to $D^k$ is determined by continuity.  The image $f_{e_1} (\partial D^k)$ is in the $(k-1)$-skeleton of $Y$, because $\overline{f}_{e_1} (\partial D^k)$ is in the $(k-1)$-skeleton of $\Lambda \backslash Y$.  Then $e_1 = f_{e_1}(D^k)$ is a $k$-cell in $Y$.  For each coset in $\Lambda / F$, let $\lambda $ be a representative, and let $e_{\lambda}$ be the $k$-cell given by $\lambda \circ f_{e_1}$.  Continuing in this way, we obtain a CW decomposition of $Y$ lifting that on $\Lambda \backslash Y$, on which $\Lambda$ acts in a cellular fashion.
\end{Pf}

This key proposition about rational cohomological dimension is known, but we did not find a reference.  

\begin{propn}
\label{folkcdq}
Let $\Lambda$ be a discrete group acting on a contractible CW complex $Y$ properly and cellularly.  Then
$$ \mbox{cd}_{\Q} \Lambda \leq \mbox{dim} Y$$
\end{propn}

\begin{Pf}

Let $C_*$ be the chain complex over $\Q$ for $Y$; this complex is exact because $Y$ is contractible.  For each dimension $n$, the $\Q \Lambda$-module $C_n$ is a direct sum
$$C_n \cong \bigoplus_{\alpha} \Q [\Lambda / \Lambda_{\alpha}] $$ 

where $\alpha$ ranges over $\Lambda$-orbits, and $\Lambda_{\alpha}$ is the stabilizer of a point in the orbit labelled by $\alpha$.  Properness of the $\Lambda$-action on $Y$ implies each $\Lambda_{\alpha}$ is finite.  There is a $\Q \Lambda$-homomorphism $\Q [ \Lambda/ \Lambda_{\alpha}]  \rightarrow \Q \Lambda$ given by
$$ [ \lambda ] \mapsto \frac{1}{ \left| \Lambda_{\alpha} \right| } \sum_{\sigma \in [ \lambda ]} \sigma$$

where $[\lambda]$ is a coset in $\Lambda / \Lambda_{\alpha}$.  This homomorphism is a section of the natural projection $\Q \Lambda \rightarrow \Q [ \Lambda / \Lambda_{\alpha}]$.  Now it follows that each $\Q [ \Lambda / \Lambda_{\alpha}]$ is a projective $\Q \Lambda$-module and thus so is each $C_n$.  Therefore, $C_*$ is a projective resolution of $\Q$ over $\Q \Lambda$.  For any $\Q \Lambda$-module $A$, the homology of the complex $Hom_{\Q \Lambda}( C_*, A)$ computes the cohomology $H^*( \Lambda, A)$.  Because $C_n = {\bf 0}$ for $n > \mbox{dim} Y$, the cohomology $H^n( \Lambda, A) = {\bf 0}$ for $n > \mbox{dim} Y$.

\end{Pf} 

\subsection{Symmetric spaces}
\label{section.symmetric.spaces}

Recall that any connected Lie group $G$ has a maximal compact subgroup $K$, unique up to conjugacy.  This subgroup is always connected; further, the quotient $G/K$ is contractible (\cite{Iwa} 6).

We collect here some facts about symmetric spaces of noncompact type, which are homogeneous Riemannian manifolds of the form $G/K$, where $G$ is semisimple with no compact local factors, connected, and has finite center.  

\begin{propn}
\label{symmspfacts}
Let $G$ be a connected semisimple Lie group with finite center.  Let $K$ be a maximal compact subgroup of $G$.  
\begin{enumerate}
\item{There is an $Ad(K)$-invariant decomposition $\mathfrak{g} = \mathfrak{p} \oplus \mathfrak{k}$.  The $Ad(K)$-irreducible subspaces of $\mathfrak{p}$ correspond to the simple factors of $\mathfrak{g}$, and there is no one-dimensional $Ad(K)$-invariant subspace of $\mathfrak{p}$.}
\item{$N_G(K) = K$}
\item{For any torsion-free lattice $\Gamma \subset G$, the center $Z(\Gamma) = 1$.}
\end{enumerate}
\end{propn}

\begin{Pf}
\begin{enumerate}
\item{The complement $\mathfrak{p}$ is the $-1$ eigenspace of a Cartan involution $s$ of $\mathfrak{g}$, and satisfies $[ \mathfrak{p}, \mathfrak{p}] \subseteq \mathfrak{k}$.  Let $B$ be the inner product on $\mathfrak{p}$ pulled back from $T_{[K]} (G/K)$.  Then $(\mathfrak{g}, s, B)$ is an orthogonal involutive Lie algebra (see \cite{Wolf} 8.2).  Suppose $\mathfrak{p}_1$ is an $Ad(K)$-invariant subspace of $\mathfrak{p}$.  Then $\mathfrak{p}_1^{\perp}$ is also $Ad(K)$-invariant.  Lemma 8.2.1 (iii) of \cite{Wolf} says that $[ \mathfrak{p}_1 , \mathfrak{p}_1^{\perp}] = {\bf 0}$.  Let $\mathfrak{k}_1$ be the centralizer in $\mathfrak{k}$ of $\mathfrak{p}_1^{\perp}$.  We will show that $\mathfrak{g}_1 = \mathfrak{p}_1 + \mathfrak{k}_1$ is an ideal of $\mathfrak{g}$.  

First, $[\mathfrak{p}_1, \mathfrak{p}_1] \subseteq \mathfrak{k}_1$:  the Jacobi identity plus the lemma above gives, for any $X,Y \in \mathfrak{p}_1$ and $Z \in \mathfrak{p}_1^{\perp}$,
$$ [[X,Y], Z] = - [[Z,X],Y] - [[Y,Z],X] = {\bf 0}$$

A similar application of the Jacobi identity shows $[\mathfrak{k}, \mathfrak{k}_1] \subseteq \mathfrak{k}_1$.  Now it is easy to verify that $\mathfrak{g}_1$ is an ideal.

Now suppose that $\mathfrak{p}_1$ is a one-dimensional $Ad(K)$-invariant subspace.  Because $K$ is compact, it has no nontrivial one-dimensional representations, so $[\mathfrak{k}, \mathfrak{p}_1] = {\bf 0}$.  Clearly $[\mathfrak{p}_1, \mathfrak{p}_1] = {\bf 0}$.  By the lemma above, $[\mathfrak{p}, \mathfrak{p}_1] = {\bf 0}$.  Therefore, $\mathfrak{p}_1 \subseteq \mathfrak{z}(\mathfrak{g})$, contradicting that $\mathfrak{g}$ is semisimple.}

\item{If $g \in N_G(K)$, then $[gK] \in G/K$ is fixed by $K$.  By \cite{Eber} 1.13.14 (4), there is a unique fixed point for $K$ in $G/K$.  Therefore, $[gK] = [K]$, so $g \in K$.}

\item{For each $z \in Z(\Gamma)$, 
$$ \gamma z \gamma^{-1} z^{-1} = 1 \ \mbox{for all} \ \gamma \in \Gamma$$

Because $\Gamma$ is Zariski dense in $G$ (see \cite{ZiETSG} 3.2.5), 
$$ g z g^{-1} z^{-1} = 1 \ \mbox{for all } \ g \in G$$

Therefore, $Z(\Gamma) \subseteq Z(G)$.  But $Z(G)$ is finite and $\Gamma$ is torsion-free. }


\end{enumerate}
\end{Pf}

\section{Examples when $Z(G^0)$ is infinite}
\label{section.examples}

\subsection{Varying metric type on orbits}
\label{subsection.varying.metric}

We now construct an example of a contractible Lorentz manifold $X$ admitting a free, properly discontinuous, cocompact, isometric action such that $\widetilde{O}(2,2k) \subseteq Isom^0(X)$, but the metric type of $\widetilde{O}(2,2k)$-orbits varies.  We believe that $Isom^0(X) \cong \widetilde{O}(2,2k)$, in which case this example illustrates the necessity of the hypothesis $|Z(G^0)| < \infty$ in the proper case of Theorem \ref{mainthm} in order to conclude that $G^0$-orbits are Riemannian. 

Let $K \cong \widetilde{O}(2k)$, the maximal compact subgroup of $\widetilde{O}(2,2k)$.  Let 
$$ \mathfrak{o}(2,2k) = \mathfrak{p}_1 \oplus \mathfrak{p}_2 \oplus \mathfrak{o}(2) \oplus \mathfrak{k}$$

be an $Ad(K)$-invariant refinement of a Cartan decomposition, with $\mbox{dim} \mathfrak{p}_i = k$ for $i=1,2$.  The bracket $[\mathfrak{p}_i,\mathfrak{p}_i] \subset \mathfrak{k}$ while $[\mathfrak{p}_i, \mathfrak{p}_j] \subseteq \mathfrak{o}(2)$ for $i,j=1,2, i \neq j$.  Let $B = B_1 \oplus B_2$ be a positive definite $Ad(K)$-invariant inner product on $\mathfrak{p}_1 \oplus \mathfrak{p}_2$ that is not invariant by the full maximal compact $Ad(O(2) \times O(2k))$ of $Ad(O(2,2k))$.  Let $\lambda$ be any nonzero element of $\mathfrak{o}(2)^* \cong \R^*$.  Then for any $c \in \R$, the inner product given by $B + c \lambda$ on $\mathfrak{p}_1 \oplus \mathfrak{p}_2 \oplus \mathfrak{o}(2)$ is $Ad(K)$-invariant.  Denote by $\widetilde{B}_c$ the corresponding left-invariant symmetric bilinear tensor on $\widetilde{O}(2,2k) / \widetilde{O}(2k)$, a line bundle over the Riemannian symmetric space $O(2,2k)/(O(2) \times O(2k))$.  

Let $X = ( \widetilde{O}(2,2k) / \widetilde{O}(2k)) \times \R$.  Denote by $\mathcal{V}$ the kernel of $TX \rightarrow T\R$.  The fibers $\mathcal{V}_{(x,t)}$ can be identified with $T_x (\widetilde{O}(2,2k) / \widetilde{O}(2k))$.  Denote by $Z$ the vector field on $X$ arising from the right action of the $1$-parameter subgroup of $\widetilde{O}(2,2k)$ generated by $\mathfrak{o}(2)$.  This vector field is invariant by the left action of $\widetilde{O}(2,2k)$.  Denote by $\mathcal{P}$ the left-invariant distribution on $X$ that coincides with $\mathfrak{p}_1 \oplus \mathfrak{p}_2$ at the the coset of the identity.  

Now define a Lorentz metric $\lambda$ on $X$ by 
\begin{eqnarray*}
\left. \lambda_{(z,t)} \right|_{\mathcal{V}} & = & \widetilde{B}_{\cos t} \\
\lambda_{(z,t)} \left( \frac{\partial}{\partial t}, \frac{\partial}{\partial t} \right) & = & - \cos t \\
\lambda_{(z,t)} \left( \frac{\partial}{\partial t}, Z \right) & = & \sin t \\
\mathcal{P} \perp \frac{\partial}{\partial t} & & 
\end{eqnarray*}

The metric on the fiber over $t$ is Lorentzian, degenerate, or Riemannian, according as $\cos t < 0$, $\cos t = 0$, or $\cos t > 0$.  The group $\widetilde{O}(2,2k) \subseteq Isom(X)$.  For $\Gamma_0$ a cocompact lattice in $\widetilde{O}(2,2k)$, the product $\Gamma_0 \times \Z$ acts properly discontinuously and cocompactly on $X$.

\subsection{Essential orbibundle}

In \cite{FW} 5.1, Farb and Weinberger construct a group $\Xi$ that acts smoothly, properly discontinuously, and cocompactly on $\R^n$, for which every finite-index subgroup contains a subgroup $F \cong \Z / p \Z$.  The action of such an $F$ on $\R^n$ necessarily has fixed points.  The quotient good orbifold $\Gamma \backslash \R^n$ is not finitely covered by any manifold.  They also find a class $\xi \in H^2(\Gamma, \Z)$ that vanishes in $H^2( \Gamma, \R)$ and is nontrivial on some $F$.  Then the extension
$$ \Z \rightarrow \widetilde{\Xi} \rightarrow \Xi$$

given by $\xi$ acts properly discontinuously, freely, and cocompactly on $X = \R \times \R^n$.  It is possible to choose a Riemannian metric $\mu$ on $\R^n$ such that the product metric $dt^2 + \mu$ on $X$ is $\widetilde{\Xi}$-invariant and has $Isom^0(X) \cong \R$.  The quotient $ M = \widetilde{\Xi} \backslash X$ is an orbibundle
$$ S^1 \rightarrow M \rightarrow \Xi \backslash X$$

that is not finitely covered by any smooth fiber bundle.

It is simple to make this example Lorentzian: take the product metric $ - dt^2 + \mu$.  We will next construct a contractible Lorentz manifold $X$ covering a compact essential orbibundle with $Isom^0(X)$ semisimple.  

Let $\Lambda$ be a cocompact lattice in $SU(1,k)$, and let $\widetilde{\Lambda}$ be the lift of $\Lambda$ to 
$\widetilde{SU}(1,k)$; it is also a cocompact lattice, containing $\Z \subset Z(\widetilde{SU}(1,k))$.  The group $\widetilde{\Xi}$ constructed above acts on the space $( \widetilde{SU}(1,k) / \widetilde{SU}(k)) \times \R^n$ because $\R$ acts on $\widetilde{SU}(1,k) / \widetilde{SU}(k)$ by flowing along the left-$\widetilde{SU}(1,k)$-invariant vector field $Z$ from the previous example.  In fact, the cocycle $\xi$ above determines an extension
$$ \widetilde{\Lambda} \rightarrow \Gamma \rightarrow \Xi$$

that acts isometrically on $X = \widetilde{AdS}^{2k+1} \times (\R^n, \mu)$.   Now $\Gamma$ acts freely, properly discontinuously, and cocompactly on $X$ by isometries, and the quotient manifold $M = \Gamma \backslash X$ is an orbibundle
$$ \widetilde{\Lambda} \backslash \widetilde{AdS}^{2k+1}  \rightarrow M \rightarrow \Xi \backslash \R^n$$

As in Farb and Weinberger's example, no finite cover of this orbibundle is a fiber bundle.  

\begin{propn}
For $X$ as above, $Isom^0(X) \cong \widetilde{O}^0(2,2k)$
\end{propn}

\begin{Pf}
Let $\mathcal{P}$ and $\mathcal{Q}$ be the distributions on $X$ tangent to $\widetilde{AdS}^{2k+1}$-fibers and $\R^n$-fibers, respectively.  Denote by $\sigma$ the tensor
$$ \sigma({\bf u},{\bf v}) = \langle R({\bf u},{\bf v}){\bf u}, {\bf v} \rangle$$

where $R$ is the curvature tensor of $X$.
The sectional curvature of a nondegenerate plane spanned by ${\bf u},{\bf v}$ is 
$$ S({\bf u},{\bf v}) = \frac{\sigma({\bf u},{\bf v})}{| {\bf u} \wedge {\bf v} |}$$

where $|{\bf u} \wedge {\bf v}| = \langle {\bf u},{\bf u} \rangle \cdot \langle {\bf v},{\bf v} \rangle - \langle {\bf u}, {\bf v} \rangle^2$.  
Because $\widetilde{AdS}^{2k+1}$ has constant curvature, $\sigma$ vanishes whenever ${\bf u}$ and ${\bf v}$ belong to $\mathcal{P}$ and span a degenerate plane (see \cite{ONeill} 8.28).   

Suppose that $\varphi^t$ is a one-parameter group of isometries of $X$.  Let ${\bf u}, {\bf v}$ be tangent vectors belonging to $\mathcal{P}_p$, where $p$ is a point of $X$, with ${\bf u}$ null, $\langle {\bf v}, {\bf v} \rangle = 1$, and ${\bf u} \perp {\bf v}$.  Let 
\begin{eqnarray*}
\varphi^t_*({\bf u}) & = & {\bf u}'_t + {\bf x}_t \\
\varphi^t_*({\bf v}) & = & {\bf v}'_t + {\bf y}_t 
\end{eqnarray*}

with ${\bf u}'_t, {\bf v}'_t \in \mathcal{P}$ and ${\bf x}_t, {\bf y}_t \in \mathcal{Q}$.  Suppose $span \{ {\bf x}_t, {\bf y}_t \}$ is two-dimensional for small nonzero $t$.  Then ${\bf u}'_t$ is timelike; further,
$$ \lim_{t \rightarrow 0} \frac{| {\bf x}_t \wedge {\bf y}_t |}{| {\bf u}'_t \wedge {\bf v}'_t |} = 0$$

Now 
\begin{eqnarray*}
0 & = & \sigma({\bf u},{\bf v}) \\
  & = & \sigma(\varphi^t_*({\bf u}), \varphi^t_*({\bf v})) \\
  & = & \sigma({\bf u}'_t,{\bf v}'_t) + \sigma({\bf x}_t, {\bf y}_t)
\end{eqnarray*}

Because $S({\bf u}'_t,{\bf v}'_t) = -1$, 
$$ \sigma({\bf x}_t, {\bf y}_t) = | {\bf u}'_t \wedge {\bf v}'_t |$$

so 
$$ S({\bf x}_t, {\bf y}_t) = \frac{| {\bf u}'_t \wedge {\bf v}'_t|}{| {\bf x}_t \wedge {\bf y}_t |}$$

which is unbounded as $t$ approaches $0$.  Because the metric on $\R^n$ is Riemannian, the sectional curvatures along $\mathcal{Q}$ are bounded on compact sets, and we have a contradiction.  Therefore, $span \{ {\bf x}_t, {\bf y}_t \}$ is at most one-dimensional for all $t$.  

Now suppose ${\bf x}_{\epsilon} \neq {\bf 0}$ for some $\epsilon >0 $.  Then ${\bf u}'_{\epsilon}$ is timelike, and 
$$\sigma({\bf u},{\bf v}) = \sigma({\bf u}'_{\epsilon} + {\bf x}_{\epsilon}, {\bf v}'_\epsilon + {\bf y}_\epsilon ) = \sigma({\bf u}'_\epsilon,{\bf v}'_\epsilon) $$

This is impossible, because the left side is zero, while the right is not.  Therefore, all null vectors in $\mathcal{P}$ have image again tangent to $\mathcal{P}$.  Because $\mathcal{P}$ is spanned by null vectors, $\varphi_t$ preserves $\mathcal{P}$, and therefore also $\mathcal{Q} = \mathcal{P}^\perp$.  Then $\varphi_t$ induces an isometry of $(\R^n, \mu)$, which must be trivial because $Isom(\R^n, \mu)$ is discrete.  It follows easily that $\varphi_t$ belongs to $\widetilde{O}(2,2k) = Isom^0(\widetilde{AdS}^{2k+1})$.
\end{Pf}

\begin{remark}  It would be interesting to know whether there can be a real-analytic example of an essential orbibundle.
\end{remark}

\section{Lorentz dynamics}

\subsection{Kowalsky's argument}
\label{kowalsky.argument}

In \cite{Kowthesis}, Kowalsky relates the dynamics of Lorentz-isometric actions of a semisimple Lie group $G$ with the adjoint representation on $Sym^2(\mathfrak{g}^*)$.

For each $x \in X$, there are linear maps 
\begin{eqnarray*}
f_x & : & \mathfrak{g} \rightarrow T_x X \\
f_x & : & Y \mapsto \left. \frac{\partial}{\partial t} \right|_0 e^{tY}x
\end{eqnarray*}

Differentiating $ g e^{tY}x = g e^{tY} g^{-1}(gx)$ gives the relation
$$ g_{*x} f_x(Y) = f_{gx} \circ \mbox{Ad}(g) (Y)$$

Let $<,>_x$ denote the inner product on $\mathfrak{g}$ obtained by pulling back the Lorentz inner product on $T_xX$ by $f_x$.  Since the $G^0$-action is isometric, 
$$ <Y,Z>_{gx} = <\mbox{Ad}(g^{-1}) Y, \mbox{Ad}(g^{-1}) Z>_x$$

In Kowalsky's argument, the dynamics of the nonproper group action imply that many root spaces of $\mathfrak{g}$ belong to the same maximal isotropic subspace for some $<,>_x$.  We adapt this argument to obtain the following result.  Recall that $\pi_i$ is the projection of $G$ or $\mathfrak{g}$ on the $i^{th}$ (local) factor.

\begin{propn} 
\label{kowarg}
Let $G$ be a connected semisimple group acting isometrically on a Lorentz manifold.  Suppose that for $y \in X$, there is a sequence $g_n \in G(y)$ with $\mbox{Ad}(g_n) \rightarrow \infty$.  Then $\mathfrak{g}$ has a root system $\Delta$ and an $\R$-split element $A$ such that
$$ \bigoplus_{\alpha \in \Delta, \alpha(A) > 0} \mathfrak{g}_{\alpha}$$

is an isotropic subspace for $<,>_y$.  

Suppose further that $G^0$ preserves an isotropic vector field $S^*$ along the orbit $G^0y$, and let $S \in \mathfrak{g}$ be such that $f_y(S) = S^*(y)$.  Then, with respect to $<,>_y$,  
$$ \left( \bigoplus_{\alpha(A) > 0} \mathfrak{g}_{\alpha} \right) \perp S$$
\end{propn}

\begin{Pf}
Let $g_n = \widehat{k}_n \widehat{a}_n \widehat{l}_n$ be the $KTK$ decomposition of $g_n$, where $T$ is a maximal $\R$-split torus in $G$, and $K = Ad^{-1}(Ad(K))$, for $Ad(K)$ a maximal compact subgroup of $Ad(G)$.  Let $Ad(g_n) = k_n a_n l_n$ be the corresponding decomposition in $Ad(G)$.  The condition $Ad(g_n) \rightarrow \infty$ implies $a_n \rightarrow \infty$.  Let $A_n = \ln a_n$.  By passing to a subsequence, we may assume
\begin{itemize}
\item{ $A_n/|A_n| \rightarrow A $ for some $\R$-split $A \in \mathfrak{g}$}
\item{$k_n \rightarrow k$}
\item{$l_n \rightarrow l$}
\end{itemize}

Let $\Delta$ be a root system with respect to $\mathfrak{a} = \ln T$.  Let $\alpha, \beta \in \Delta$ be such that $\alpha(A), \beta(A) > 0$.  Let $U \in \mathfrak{g}_{\alpha}$ and $V \in \mathfrak{g}_{\beta}$.  We have, for all $n$,
\begin{eqnarray}
< U,V>_{\hat{a}_n \hat{l}_n y} = <U,V>_{\hat{k}_n^{-1} y}  \label{roots.iso} 
\end{eqnarray}

The left hand side is
\begin{eqnarray*}
< a_n^{-1} (U), a_n^{-1} (V)>_{\hat{l}_n y} & = & e^{- \alpha(A_n) - \beta(A_n)} <U,V>_{\hat{l}_n y} \\
  & = & e^{- \alpha(A_n) - \beta(A_n)} < l_n^{-1}(U), l_n^{-1}(V) >_y \\
\end{eqnarray*}

The inner products $< l_n^{-1}(U), l_n^{-1}(V) >_y$ converge to $<l^{-1}(U), l^{-1}(V)>_y$; in particular, they are bounded.  The factors $e^{ - \alpha(A_n) - \beta(A_n)}$ converge to $0$.  Then the left side of (\ref{roots.iso}) converges to $0$. 

The right hand side of (\ref{roots.iso}) converges to 
$$ < k (U), k (V)>_{y}$$

Therefore, the sum of root spaces
$$ \bigoplus_{\alpha(A) > 0} k(\mathfrak{g}_{\alpha})$$

is an isotropic subspace for $<,>_y$.  Now replace $\Delta$ with $\Delta \circ k^{-1}$ and $A$ with $ k (A)$ to obtain the first assertion of the proposition. 

Now let $S^*$ be a $G^0$-invariant vector field along the orbit $G^0y$.  If $S$ is such that $f_y (S) = S^*(y)$, then $f_{gy}(Ad(g)(S)) = S^*(gy)$ for any $g \in G^0$.  Let $Ad(g_n) = k_n a_n l_n$ be the $KTK$ decomposition as above, and let $A = \lim(A_n / |A_n|)$.  Now suppose $\alpha$ is a root with $\alpha(A) > 0$.  For $U \in \mathfrak{g}_{\alpha}$
$$ < k_n U, k_n a_n l_n (S)>_{\widehat{k}_n \widehat{a}_n \widehat{l}_n y} = < k_n U, S>_{y}$$

The left hand side is 
$$ < a_n^{-1}(U), l_n (S)>_{\widehat{l}_n y}  =  e^{- \alpha(A_n)} < U, l_n(S)>_{\widehat{l}_n y}$$

This sequence converges to $0$.  The right hand side converges to 
$$ < k (U), S>_y $$

Then $k (\mathfrak{g}_{\alpha}) \perp S$ with respect to $<,>_y$, yielding the desired result when $A$ is replaced with $k (A)$ and $\Delta$ with $\Delta \circ k^{-1}$.
\end{Pf}

\begin{remark}
\label{kowarg.projn.remark}
Note that if, for the $\R$-split element $A$ given by Proposition \ref{kowarg}, $\pi_i(A) \neq {\bf 0}$, then $\pi_i(g_n) \rightarrow \infty$.
\end{remark}

\begin{remark} 
\label{kak.remark}
In the proof above, if we start with a $KTK$ decomposition with $\mathfrak{a} = \ln T$, then the element $A$ given by Proposition \ref{kowarg} belongs to $Ad(K)(\mathfrak{a})$.
\end{remark}

\subsection{Totally geodesic codimension-one lightlike foliations}

A \emph{lightlike} submanifold of a Lorentz manifold is a submanifold on which the restriction of the metric is degenerate.  A foliation is lightlike if each leaf is lightlike.  In \cite{Zetgl1}, Zeghib shows that a compact Lorentz manifold $M$ with a noncompact group $G \subset Isom(M)$ has totally geodesic codimension-one lightlike (tgl) foliations.  Fix a smooth Riemannian metric $\sigma$ on $M$ giving rise to a norm $| \cdot |$ and a distance $d$ on $M$.  Let $x \in M$ and $g_n$ be a sequence in $G$.  The \emph{approximately stable set of $g_n$ at $x$} is
$$AS(x,g_n) = \{ v \in T_xM : v = \lim v_n \ \mbox{where} \ v_n \in TM \ \mbox{and}\ |g_{n*} v_n| \ \mbox{is bounded} \}$$

Zeghib proves that any unbounded $g_n$ has a subsequence for which the approximately stable set in $TM$ forms an integrable codimension-one lightlike distribution with totally geodesic leaves.  The resulting foliation $\mathcal{F}$ is Lipschitz, in the sense that there exists $C > 0$ such that
$$\angle(T \mathcal{F}_x, T \mathcal{F}_y) \leq C \cdot d(x,y)$$

for all sufficiently close $x,y \in M$.  Provided $x$ and $y$ are in a common normal neighborhood, we can define the angle above as 
$$\angle(T \mathcal{F}_x, T \mathcal{F}_y ) = \angle_{\sigma} ( P_{\gamma} T_x \mathcal{F}_x, T_y \mathcal{F}_y )$$

where $P_{\gamma}$ is parallel transport with respect to the Lorentzian connection along the geodesic $\gamma$ from $x$ to $y$.  In fact, there exists $C$ that serves as a uniform Lipschitz constant for all totally geodesic codimension-one foliations.

We extend this work to obtain tgl foliations on $X$ associated to a sequence $g_n \in G$ unbounded modulo $\Gamma$.  Let $|\cdot |$ be a smooth norm on $X$ that is $\Gamma$-invariant; such a norm can be obtained by lifting an arbitrary smooth norm from $M$.  For $x \in X$ and a sequence $g_n \in G$, define 
$$AS(x,g_n) = \{ v \in T_xX : v = \lim v_n \ \mbox{where} \ v_n \in TX \ \mbox{and}\ |g_{n*} v_n| \ \mbox{is bounded} \}$$

Note that $AS(x,g_n) = AS(x,\gamma_n g_n)$ for any sequence $\gamma_n$ in $\Gamma$, so this set can be considered associated to a sequence in $\Gamma \backslash G$.  On the other hand, for $g \in G$, 
$$AS(x,g_n g^{-1}) = g_* (AS(x,g_n))$$

\begin{thm}
\label{tglsexist}
Let $g_n \in G$ be unbounded modulo $\Gamma$.  Then there is a subsequence such that the set of $AS(x,g_n)$, for $x \in X$, form an integrable distribution with totally geodesic codimension-one lightlike leaves.  Moreover, the set $\mathcal{TGL}(X)$ of tgl foliations is uniformly Lipschitz: there exist $C, \delta > 0$, such that, for any foliation $\mathcal{F} \in \mathcal{TGL}(X)$, for any $x,y \in X$ with $d(x,y) < \delta$,
$$ \angle( \mathcal{F}_x, \mathcal{F}_y) \leq C \cdot d(x,y) $$
 \end{thm}

 The proof is essentially the same as that in \cite{Zetgl1}.  We outline that proof here and provide the observations relevant to our generalization.  For completeness, we prove the uniformly Lipschitz property in detail in the Appendix.  

  We begin with a definition.

\begin{defn}
\label{defn.geod.lamn}
Let $X$ be a $k$-dimensional manifold endowed with a smooth, torsion-free connection $\nabla$ and a smooth Riemannian metric $\sigma$.  A \emph{radius-$r$ codimension-one geodesic lamination} on $X$ consists of a subset $X^{\prime} \subset X$ and a section $f : X^{\prime} \rightarrow Gr^{k-1}(TX)|_{X^{\prime}}$, satisfying
\begin{enumerate}
\item{$\mathcal{L}_x = exp^{\nabla}(f(x) \cap B_{\sigma}({\bf 0},r))$ is $\nabla$-geodesic for each $x \in X^{\prime}$}
\item{$\mathcal{L}_x \cap \mathcal{L}_y$ is open in both $\mathcal{L}_x$ and $\mathcal{L}_y$ for all $x, y \in X^{\prime}$}
\end{enumerate}
\end{defn}

\begin{propn}
\label{uniformly.lipschitz}
Let $X$ be the universal cover of a compact manifold $M$.  Let $\nabla$ be a smooth connection and $\sigma$ a smooth Riemannian metric, both lifted from $M$.   For any $r>0$, there exist $C, \delta > 0$ such that any radius-$r$, codimension-one geodesic lamination $(X^{\prime},f)$ on $X$ is $(C,\delta)$-Lipschitz:  any $x,y \in X^{\prime}$ with $d_{\sigma}(x,y) < \delta$ are connected by a unique $\nabla$-geodesic $\gamma$, and 
$$ \angle_{\sigma} (P_{\gamma} f(x), f(y)) \leq C \cdot d_{\sigma}(x,y)$$
\end{propn}

The proof of this proposition is in the Appendix.  We record two consequences.

\begin{cor}
\label{geod.lamn.unif.cts}
For any radius-$r$ codimension-one geodesic lamination $(X^{\prime}, f)$, the function $f$ is uniformly continuous on $X^{\prime}$.
\end{cor}

\begin{cor}
\label{geod.lamn.compact}
The space $\mathcal{TGL}(X)$ is compact.
\end{cor}

For the remainder of this section, all metric notions, such as the distance $d$ and norm $| \cdot |$, always refer to $\sigma$ below.  Affine notions, such as geodesics, parallel transport $P$, and the exponential map $exp$, refer to $\nabla$.

\begin{Pf} (of Theorem \ref{tglsexist})  
 All exact references are to \cite{Zetgl1}, unless otherwise indicated.

The following terminology will be used below.  A sequence of subspaces $H_n \subset TX$ is a \emph{stable sequence} for $g_n \in G$ if $ \parallel \left. g_{n*} \right|_{H_n} \parallel$ is bounded.  The \emph{modulus of stability} of a stable sequence $H_n$ is $1/(sup_n  \parallel \left. g_{n*} \right| _{H_n} \parallel) $.

For $x \in X$ and $g_n \in G$, the \emph{punctually approximately stable set of $g_n$ at $x$} is
$$ PAS(x,g_n) =  \{v \in T_xX : v = \lim v_n \ \mbox{where} \ v_n \in T_xX \ \mbox{and}\ |g_{n*} v_n| \ \mbox{is bounded} \}$$

{\bf Step 1}: Punctually approximately stable hyperplanes on a dense subset.
\nopagebreak[3]

Endow $\R^k$ with the standard Lorentz inner product and the standard positive-definite quadratic form $Q$.  Let $\overline{\theta} : TM \rightarrow M \times \R^k$ be a measurable, almost-everywhere smooth, bounded trivialization of $TM$.  Bounded in this case means there exists $c > 0$ such that 
$$ |v| /c \leq Q(\overline{\theta}(x, v))^{1/2} \leq c |v| $$

for any $x \in M$ and $v \in T_xM$.  Such a $\overline{\theta}$ can be obtained from finitely many local trivializations covering $TM$ so that $\overline{\theta}$ is smooth on the complement of finitely many spheres in $M$.

Lift $\overline{\theta}$ to a $\Gamma$-equivariant trivialization $\theta : TX \rightarrow X \times \R^k$, where $\Gamma$ acts on $X \times \R^k$ by $\gamma .(x, {\bf v}) = (\gamma.x, {\bf v})$.  This trivialization of $TX$ is also bounded.  Denote by $\alpha$ the resulting cocycle
$$ \alpha : G \times X \rightarrow O(1,k-1)$$

Denote by $\theta_x$ the restriction of $\theta$ to $T_xX$.  It is an isomorphism $T_xX \rightarrow \R^k$.

Let $g_n$ be a sequence in $G$ and $x \in X$ such that $\alpha(g_n,x)$ is defined for all $n$.  If $\parallel (g_n)_{*x} \parallel \rightarrow \infty$, then $A_n = \alpha(g_n, x) \rightarrow \infty$ in $O(1,k-1)$.  Corollary 4.3 states that any $A_n \rightarrow \infty$ in $O(1,k-1)$ has a subsequence $B_n$ for which $AS({\bf 0}, B_n)$ is a lightlike hyperplane $H$.  For such a sequence $B_n$, the stable sequence $H_n$ of hyperplanes can be chosen so that the modulus of stability is $1$.  After passing to the corresponding subsequence of $g_n$, the punctually approximately stable space is
$$PAS(x,g_n) = (\theta_x)^{-1}(H) = \lim (\theta_x)^{-1}(H_n)$$

Because $\theta$ is bounded, the modulus of stability is uniformly bounded for all such sequences.  Let this bound be $r > 0$.

Let $U$ be the open dense subset on which $\theta$ is defined.  For any sequence $g_n$, the intersection $D = (\cap_n g_nU) \cap U$ has conull measure in $X$.  Let $X^{\prime} \subset D$ be a countable dense subset of $X$.  For each $x \in X^{\prime}$, the cocycle $\alpha(g_n,x)$ is defined for all $n$.  There is a subsequence such that, for all $x \in X^{\prime}$, the punctually approximately stable sequence $PAS(x,g_n)$ is a lightlike hyperplane.  

{\bf Step 2}: Propagation of stability.
\nopagebreak[3]

On $M$, there exists $r^{\prime} > 0$ such that the exponential map is defined on $B^{r^{\prime}}M \subseteq TM$.  By $\Gamma$-invariance of $\sigma$ and of the Lorentz metric, the exponential map for $X$ is also defined on $B^{r^{\prime}}X$.  We may assume that our uniform modulus of stability $r \leq r^{\prime}$.  

Proposition 6.1 says that if $H_n$ is a stable sequence of linear spaces for $g_n$, then on the submanifolds $\mathcal{L}_n = exp(H_n \cap B^r X)$, the derivatives $\left. g_{n*} \right| _{T \mathcal{L}_n} $ are uniformly bounded.  The key fact is that if $M$ is compact, then the derivatives of 
$$ exp \circ \overline{\theta}^{-1} : M \times \R^k \rightarrow M \times M$$

are uniformly bounded.  In our case, because $exp$ and $\theta$ are $\Gamma$-equivariant, and $\Gamma$ preserves $\sigma$, the same is true.

Proposition 6.2 says that for $H_n$ as above, with $H_n \subseteq T_{x_n} X$, for any $C^1$-bounded sequence of curves $c_n : [0,1] \rightarrow \mathcal{L}_n$ with $c_n (0) = x_n$, the parallel transports $P_{c_n} H_n$, also form a stable sequence for $g_n$.  This proposition is local.  The proof uses the previous proposition and the fact that $g_{n*}$ commutes with parallel transport.

{\bf Step 3}: Geodesibility.
\nopagebreak[3]

Proposition 6.3 says that if $ \lim H_n = H \subset T_xM$ is an approximately stable subspace, not properly contained in any other approximately stable subspace, then $\mathcal{L}_x = exp(H \cap B^rX)$ is totally geodesic.  This proposition is local. The proof uses the propositions of Step 2.  

{\bf Step 4}: Approximately stable subspaces have positive codimension.
\nopagebreak[3]

Corollary 5.2 says that for a compact Lorentz manifold $M$, if $g_n \in Isom(M)$ has $AS(x,g_n) = T_xM$ for some $x \in M$, then $g_n$ is bounded.  Using Proposition 6.1 of Step 2, one can see that $AS(x,g_n) = T_xM$ implies there is a stable sequence $T_{x_n} M$ for $g_n$.  Then the derivatives $(g_n)_{* x_n}$ are bounded, which implies $g_n$ is equicontinuous if $M$ is compact.  In fact, for an arbitrary Lorentz manifold $X$, if $K \subset X$ is compact and there exist $x_n \in K$ such that $g_n x_n \in K$ and $(g_n)_{*x_n}$ are bounded, then $g_n$ is equicontinuous.

Now let $X$ be the universal cover of a compact Lorentz manifold $M = \Gamma \backslash X$.  Suppose $g_n$ is a sequence with $AS(x,g_n) = T_xX$ for some $x$, so there is a stable sequence $T_{x_n} X$ for $g_n$.  Let $K$ be a compact fundamental domain for $\Gamma$ with $x_n,x \in K$, passing to a subsequence if necessary.  Let $\gamma_n \in \Gamma$ be such that $\gamma_n g_n x_n \in K$.  Then $AS(x, \gamma_n g_n ) = AS(x, g_n)$, so $\gamma_n g_n$ is bounded.  We conclude that, if $g_n$ is unbounded modulo $\Gamma$, then, for all $x \in X$, the approximately stable set $AS(x,g_n)$ is a proper subset of $T_xX$. 

{\bf Step 5}: Totally geodesic codimension-one laminations.
\nopagebreak[3]

Fact 6.4 says that if $g_n$ does not admit a codimension-zero approximately stable subspace and $H \subseteq T_xX$ is an approximately stable hyperplane, then $AS(x,g_n) = H$.  The proof is local and uses Proposition 6.1 of Step 2.  From Steps 1 and 4, we have that if $g_n$ is unbounded modulo $\Gamma$, then there is a subsequence $g_n$ for which $PAS(x,g_n) = AS(x,g_n)$ is a lightlike hyperplane.  In fact, for a countable dense subset $X^{\prime} \subset X$, there is a subsequence $g_n$ such that $PAS(x,g_n) = AS(x,g_n)$ is a lightlike hyperplane for all $x \in X^{\prime}$.  From Proposition 6.3 of Step 3, $\mathcal{L}_x = exp(AS(x,g_n) \cap B^r X)$ is a totally geodesic hypersurface for all $x \in X^{\prime}$.

Corollary 6.5 says that, for any $y \in \mathcal{L}_x$, the hyperplane $T_y \mathcal{L}_x = AS(y, g_n)$.  Corollary 6.6 says that if $\mathcal{L}_x$ and $\mathcal{L}_y$ are two totally geodesic hypersurfaces arising in this way, then $\mathcal{L}_x \cap \mathcal{L}_y$ is open in both $\mathcal{L}_x$ and $\mathcal{L}_y$.  Both follow from Proposition 6.1 of Step 2 and Step 4. 

Now for any $g_n$ unbounded modulo $\Gamma$, there is a subsequence $g_n$ such that $(X^{\prime}, f)$, where $f(x) = PAS(x,g_n)$, is a radius-$r$ codimension-one totally geodesic lamination.

{\bf Step 6}: Totally geodesic codimension-one lightlike foliations.
\nopagebreak[3]

From our Corollary \ref{geod.lamn.unif.cts} above, the map $f : x \mapsto PAS(x,g_n)$ is uniformly continuous on $X^{\prime}$.  Given $x \in X$, let $x_m \in X^{\prime}$ be a sequence converging to $x$.  The limit $H = \lim_{m \rightarrow \infty} PAS (x_m,g_n)$ exists and is a lightlike hyperplane.  Because the modulus of stability for the stable sequence converging to $PAS(x_m,g_n)$ is uniform over $m$, the limit $H$ is approximately stable for $g_n$.  Then Fact 6.4 of Step 5 gives that $H = AS(x,g_n)$.  Therefore, the lamination $(X^{\prime},f)$ extends to a totally geodesic codimension-one lightlike foliation $\mathcal{F}$ with $\mathcal{F}_x = AS(x,g_n)$ for all $x \in X$.  \end{Pf}

\section{A closed orbit}

In this section, we consider a slightly more general setting. 
Let $M$ be a compact, connected, real-analytic manifold with a real-analytic rigid geometric structure of algebraic type defining a connection.  Assume that this connection is \emph{complete}---that is, that the exponential map is defined on all of $TM$.  An example of a rigid geometric structure of algebraic type defining a connection is a pseudo-Riemannian metric.   We will make use of Gromov's stratification theorem and its consequences for real-analytic rigid geometric structures of algebraic type.  Let $G$ be the group of automorphisms of the lifted structure on the universal cover $X$; it is a finite-dimensional Lie group (\cite{Gromov} 1.6.H).  As usual, let $G^0$ be the identity component of $G$; let $\Gamma \subset G$ be the group of deck transformations of $X$; and let $\Gamma_0 = \Gamma \cap G^0$.  

  Let $J$ be the pseudogroup of germs of local automorphisms of $M$.  For $x \in M$, let $J_x$ be the pseudogroup of germs at $x$ of local isometries.  Call the $J$-\emph{orbit} of $x \in M$ the equivalence class of $x$ under the relation $x \sim y$ when $jx = y$ for some $j \in J_x$.  Gromov's stratification theorem says the following:

\begin{thm}[\cite{Gromov} 3.4]   There is a $J$-invariant stratification
\label{stratification}
$$ \emptyset = M_{-1} \subset M_0 \subset \cdots \subset M_k = M$$

such that, for each $i$, $0 \leq i \leq k$, the complement $M_i \backslash M_{i-1}$ is an analytic subset of $M_i$.  Further, each $M_i \backslash M_{i-1}$ is foliated by $J$-orbits, and the $J$-orbits are properly embedded in $M_i \backslash M_{i-1}$.
\end{thm}

\begin{cor}[\cite{Gromov} 3.4.B, compare \cite{DAG} 3.2.A (iii)] 
\label{closed.orbit}
There exists a closed $J$-orbit in $M$.
\end{cor}


The stratification above is obtained from similar stratifications invariant by infinitesimal isometries of order $k$, for arbitrary sufficiently large $k$.  It is shown in \cite{Gromov} 1.7.B that orbits of infinitesimal isometries of increasing order eventually stabilize to $J$-orbits.  For any $x \in M$, the infinitesimal isometries of order $k$ fixing $x$ form an algebraic subgroup of $GL(T_xM)$, because the given $H$-structure is of algebraic type.  Then stabilization of infinitesimal isometries to local isometries implies that the group $J(x)$ of germs in $J_x$ fixing $x$ has algebraic isotropy representation on $T_xM$ (see \cite{DAG} 3.5, \cite{Gromov} 3.4.A); in particular, $J(x)$ has finitely-many components.

The aim of this section is to establish that the properties of $J$-orbits discussed above apply also to images in $M$ of $G^0$-orbits on $X$.  The main reason for this correspondence is the fact, proved by Nomizu \cite{Nomizu}, Amores \cite{Amores}, and, in full generality, Gromov \cite{Gromov}, that local Killing fields on $X$ can be uniquely extended to global Killing fields.  Because the connection on $X$ is complete, any global Killing field integrates to a one-parameter subgroup of $G$ (see \cite{KN} VI.2.4).  Thus there is a correspondence between local Killing fields near any point of $M$ and elements of $\mathfrak{g}$.

\begin{propn}
\label{algstabs}
The isotropy representation of $G(y)$ is algebraic for any $y \in X$; the same is true for $G^0(y)$.
\end{propn}

\begin{Pf}
Denote by $\pi$ the covering map from $X$ to $M$. There is an obvious homomorphism $\varphi: G(y) \rightarrow J(z)$, where $z = \pi(y)$.  A tangent vector at the identity to $J(z)$ corresponds to the germ of a local Killing field at $z$.  Local Killing fields near $z$ can be lifted to $X$, extended, and integrated, giving a linear homomorphism $T_e(J(z)) \rightarrow \mathfrak{g}$ inverse to $D_e \varphi$.  Then $\varphi$ is a local diffeomorphism near the identity, and so it is a local isomorphism $G(y) \rightarrow J(z)$.  By rigidity, any $g \in G(y)$ with trivial germ at $y$ is trivial, so $\varphi$ is an isomorphism onto its image.  The image is a union of components of $J(z)$, so $G(y)$ is algebraic.  The restriction of $\varphi$ to $G^0(y)$ is also an isomorphism onto its image.
\end{Pf}

\begin{propn}
\label{goodorbit}
There is an orbit $G^0y$ in $X$ with closed image in $M$.
\end{propn}

\begin{Pf}
Let $z \in M$ have closed $J$-orbit, and choose any $y \in X$ with $\pi(y) = z$.  The image $\pi(G^0y)$ is a connected submanifold of $Jz$, though it is not \emph{a priori} closed.  Denote by $J^0z$ the component of $z$ in $Jz$.  This is the orbit of $z$ under local Killing fields on $M$---that is, all points of $M$ that can be reached from $z$ by flowing along a finite sequence of local Killing fields.  Because each local Killing field on $M$ corresponds to a $1$-parameter subgroup of $G^0$, this component $J^0z$ is contained in $\pi(G^0y)$.  They are therefore equal, and closed in $M$, because $Jz$ has finitely-many components and is closed in $M$.
\end{Pf}

\begin{propn}
\label{gamma0.cocompact}
Let $y$ be as in the previous proposition, so $G^0y$ has closed image in $M$.
The subgroup $\Gamma_0 = G^0 \cap \Gamma \subset G^0$ acts freely, properly discontinuously, and cocompactly on $G^0/G^0(y)$.
\end{propn}

\begin{Pf}
Let $G_y$ be the subgroup of $G$ leaving invariant the orbit $G^0 y$, and $\Gamma_y = G_y \cap \Gamma$; note that $\Gamma_y$ acts cocompactly on $G^0 y$.  Because $G^0y$ is a closed submanifold of $X$, the orbit map $G^0/G^0(y) \rightarrow G^0 y$ is a homeomorphism onto its image.  It therefore suffices to show that $\Gamma_0$ has finite index in $\Gamma_y$.  


Now $G^0 y = G_y y$ is also the homeomorphic image of $G_y / G(y)$, which is then connected.  As in Proposition \ref{algstabs}, $G(y)$ has finitely-many components; then so does $G_y$.  Thus $G^0$ is a finite-index subgroup of $G_y$, so $\Gamma_0$ is a finite-index subgroup of $\Gamma_y$, as desired.  
\end{Pf}

\begin{cor} If $G^0$ has no compact orbits on $X$, then $\Gamma_0$ is an infinite normal subgroup of $\Gamma$.
\end{cor}

\section{Proof of main theorem}

Let $Y = G^0 y$ be the orbit given by Proposition \ref{goodorbit} with closed projection to $M$. 

\subsection{Proper case}
\label{section.proper.case}
If $G^0(y)$ is compact, then $G^0$ acts properly; in fact, so does the group $G^{\prime}$ generated by $G^0$ and $\Gamma$.

\begin{propn}
Let $G^{\prime}$ be the closed subgroup of $G$ generated by $G^0$ and $\Gamma$.  If $G^0(y)$ is compact, then $G^{\prime}$ acts properly on $X$.
\end{propn}

\begin{Pf}
If $G^0(y)$ is compact, then $\Gamma_0$ is a cocompact lattice in $G^0$ by Proposition \ref{gamma0.cocompact}.  Let $F$ be a compact fundamental domain for $\Gamma_0$ containing the identity in $G^0$; note $F$ is also a compact fundamental domain for $\Gamma$ in $G^{\prime}$.  Let $A$ be a compact subset of $X$, and $G^{\prime}_A$ the set of all $g$ in $G^{\prime}$ with $gA \cap A \neq \emptyset$.  Any $g \in G^{\prime}_A$ is a product $\gamma f$ where $f \in F$ and $\gamma \in \Gamma_{FA}$.  Since $FA$ is compact, $\Gamma_{FA}$ is a finite set $\{ \gamma_1, \ldots, \gamma_l \}$.  Then $G_A^{\prime}$ is a closed subset of the compact set $\gamma_1 F \cup \cdots \cup \gamma_l F$, so it is compact.  
\end{Pf}

The first statement in the proper case of Theorem \ref{mainthm} is that $M$ is an orbibundle
$$ \Lambda \backslash G^0 / K_0 \rightarrow M \rightarrow Q$$

We prove this statement, with $\Lambda = \Gamma_0$, in three steps.

{\bf Step 1}: $\Gamma / \Gamma_0$ proper on $G^0 \backslash X$.
\nopagebreak[3]

Let $\overline{A}$ be a compact subset of $G^0 \backslash X$, and let 
$$(\Gamma / \Gamma_0)_{\overline{A}} = \{ [ \gamma ] \in \Gamma / \Gamma_0 : [ \gamma ] \overline{A} \cap \overline{A} \neq \emptyset \}$$

The aim is to show this set is finite.
There is a compact subset $A$ of $X$ projecting onto $\overline{A}$ by Proposition \ref{stratfacts} (1).  Let 
$$ \Gamma_{A,G^0A} = \{ \gamma \in \Gamma : \gamma A \cap G^0 A \neq \emptyset \}$$  

Note that $\Gamma_{A,G^0A}$ is invariant under right multiplication by $\Gamma_0$, and
$$(\Gamma / \Gamma_0)_{\overline{A}} = \Gamma_{A,G^0A} / \Gamma_0$$

Let $F$ be a compact fundamental domain for $\Gamma_0$ in $G^0$.  Since $FA$ is compact and $\Gamma$ acts properly, the set $\Gamma_{A,FA}$ is finite.  Then $(\Gamma_{A,FA} \cdot \Gamma_0) / \Gamma_0 = \Gamma_{A,G^0A} / \Gamma_0$ is finite, as well.

\medskip

Let $K_0$ be a maximal compact subgroup of $G^0$.  

{\bf Step 2}: $G^0(x) \cong K_0$ for all $x \in X$.
\nopagebreak[3]

Any stabilizer $G^0(x)$ is compact, so conjugate to a subgroup of $K_0$.  Since $K_0$ is connected, it suffices to show that $\mbox{dim} K_0 \leq \mbox{dim} G^0(x)$.  We follow the cohomological dimension arguments of Farb and Weinberger (\cite{FW}).  By Proposition \ref{basiccdq} (2), 
$$\mbox{cd}_{\Q} \Gamma = \mbox{dim} X \qquad \mbox{and} \qquad \mbox{cd}_{\Q} \Gamma_0 = \mbox{dim} (G^0 / K_0)$$

By the extension of \cite{FW} (2.2) of the Conner conjecture (\cite{Oliver}), the quotient space $G^0 \backslash X$ is contractible because $G^0$ acts properly and $X$ is contractible. For any $x \in X$, 
$$ \mbox{dim} X - \mbox{dim} (G^0 / G^0(x)) \geq \mbox{dim}(G^0 \backslash X)$$

by \ref{stratfacts} (2).  The quotient $(\Gamma / \Gamma_0) \backslash (G^0 \backslash X) =  G^{\prime} \backslash X$ is a Whitney stratified space because $G^{\prime}$ acts properly (\ref{stratfacts} (2)), so it is triangulable by \ref{strattriang}.  By Proposition \ref{cellular}, there is a CW decomposition of $G^0 \backslash X$ preserved by the $(\Gamma / \Gamma_0)$-action,  so \ref{folkcdq} gives
$$\mbox{cd}_{\Q} (\Gamma / \Gamma_0) \leq \mbox{dim} (G^0 \backslash X) $$  

Now the inequality \ref{basiccdq} (3) gives 
\begin{eqnarray*}
\mbox{dim} X & \leq & \mbox{dim}(G^0/K_0) + \mbox{dim}(G^0 \backslash X) \\
 & \leq & - \mbox{dim}K_0 + \mbox{dim}X + \mbox{dim}G^0(x)
\end{eqnarray*}

for any $x \in X$, as desired.  

\medskip

{\bf Step 3}: Orbibundle.
\nopagebreak[3]

Now $G^0 \backslash X$ is a manifold (\ref{stratfacts} (3)) on which $\Gamma / \Gamma_0$ acts properly discontinuously.  The foliation of $X$ by $G^0$-orbits descends to $M$, and all leaves in $M$ are closed.  The leaf space is $Q = (\Gamma / \Gamma_0) \backslash (G^0 \backslash X)$, a smooth orbifold.  Given $U$ open in $Q$, lift it to a connected $\widetilde{U}$ in $G^0 \backslash X$.  For $U$ sufficiently small, the fibers of $M$ over $U$ are 
$$ \widetilde{U} \times_{\Lambda_{\widetilde{U}}} \Gamma_0 \backslash G^0 / K_0$$

where $\Lambda_{\widetilde{U}} = \{ [\gamma] \in \Gamma / \Gamma_0 \ : \ [\gamma] \widetilde{U} \cap \widetilde{U} \neq \emptyset \}$ is a finite group.  We have an orbibundle
$$ \Gamma_0 \backslash G^0 / K_0 \rightarrow M \rightarrow Q$$

Now it remains to prove the second part of the theorem in the proper case, giving the metric on $M$, assuming $Z(G^0)$ is finite.  

{\bf Step 4}: Splitting of $X$.
\nopagebreak[3]

Let 
\begin{eqnarray*}
 & \rho : X \rightarrow G^0 / K_0 & \\
 \rho(x) = [g] & \qquad \mbox{where} \qquad & gK_0 g^{-1} = G^0(x)
\end{eqnarray*}

This map is well-defined and injective along each orbit because $N(K_0) = K_0$ (\ref{symmspfacts} (2)).  Each fiber $\rho^{-1}([g]) = \mathit{Fix}(gK_0g^{-1})$.  Each orbit is mapped surjectively onto $G^0/ K_0$.  Let $L = \rho^{-1}([e]) = Fix(K_0)$, a totally geodesic submanifold of $X$.  Under the quotient, $L$ maps diffeomorphically to $G^0 \backslash X$, so $L$ is connected.  The map 
$$G_0 / K_0 \times L \rightarrow X \qquad ([g],l) \rightarrow g l$$

is a well-defined diffeomorphism.

The restriction of the metric to each $G^0$-orbit must be Riemannian.  Indeed, let $x \in L$ and consider the isotropy representation of $K_0$.  The map $f_x : \mathfrak{g} \rightarrow T_x X$ gives a $K_0$-equivariant isomorphism $\mathfrak{g} / \mathfrak{k} \rightarrow T_x(G^0 x)$.  If the inner-product on $T_x(G^0 x)$ is degenerate, then $K_0$ is trivial on the kernel subspace. If it is Lorentzian, then $K_0$ preserves a norm, so it fixes a minimal length timelike vector.  Either way, the isotropy representation of $K_0$ has a fixed vector.  But $\mbox{Ad}(K_0)$ has no one-dimensional invariant subspace in $\mathfrak{p} \cong \mathfrak{g} / \mathfrak{k}$ (\ref{symmspfacts} (1)), a contradiction.

For the same reason, $L$ is orthogonal to each $G^0$-orbit.  Let $x \in L$.  The subspaces $T_x(G^0 x)^{\perp}$ and $T_xL$ are both $K_0$-invariant complements to $T_x(G^0 x)$ in $T_xX$.  If they are unequal, then there are nonzero vectors $v \in T_x(G^0 x)$ and $w \in T_xL$ such that $v-w \in T_x(G^0 x)^{\perp}$. Then 
\begin{eqnarray*}
& & k (v-w) = k v - w \in T_x(G^0 x)^{\perp} \\
& \Rightarrow & k v - v  \in T_x(G^0 x)^{\perp} \\
& \Rightarrow & k v = v
\end{eqnarray*}

again contradicting that $K_0$ has no one-dimensional invariant subspace in $T_x(G^0 x)$.

{\bf Step 5}: Splitting of $\Gamma$.
\nopagebreak[3]

The argument here is the same as in \cite{FW}.  The extension
$$ \Gamma_0 \rightarrow \Gamma \rightarrow \Gamma / \Gamma_0$$

is a subextension of 
$$ G^0 \rightarrow G^{\prime} \rightarrow \Gamma / \Gamma_0$$

so the action $\Gamma / \Gamma_0 \rightarrow Out(\Gamma_0)$ is the restriction of $\Gamma / \Gamma_0 \rightarrow Out(G^0)$.  Since $G^0$ is semisimple, $Out(G^0)$ is finite.  Thus there is a finite-index subgroup $\Gamma^{\prime}$ of $\Gamma$ containing $\Gamma_0$ such that conjugation by any $\gamma \in \Gamma^{\prime}$ is an inner automorphism of $\Gamma_0$.  The extension 
$$\Gamma_0 \rightarrow \Gamma^{\prime} \rightarrow \Gamma^{\prime} / \Gamma_0$$

also determines a cocycle in $H^2(\Gamma^{\prime} / \Gamma_0, Z(\Gamma_0))$.  But $Z(\Gamma_0)$ is trivial (\ref{symmspfacts} (3)).  This extension is therefore a product 
$$\Gamma^{\prime} \cong \Gamma_0 \times \Gamma^{\prime} / \Gamma_0 $$

Since $\Gamma^{\prime}$ is torsion-free (\ref{basiccdz} (3)), so is $\Gamma^{\prime} / \Gamma_0$.  Then $\Gamma^{\prime} / \Gamma_0$ acts freely on $G^0 \backslash X$, and the quotient, which is a finite cover of $Q$, is a manifold $Q^{\prime}$.  The finite cover $M^{\prime} = \Gamma^{\prime} \backslash X$ is diffeomorphic to $\Gamma_0 \backslash G^0 / K_0 \times Q^{\prime}$.  The metric descends from $X$ to $M^{\prime}$ and has the form claimed in the theorem.
 
\subsection{Nonproper case: if $G^0$ has infinite orbit in $\mathcal{TGL}(X)$}

Now suppose that $G^0(y)$ is noncompact, so $G^0$ acts nonproperly; further, $\Gamma \backslash G$ is noncompact.  By Theorem \ref{tglsexist}, there are tgl foliations on $X$.  The set $\mathcal{TGL}(X)$ of all these foliations forms a $G$-space.  Pick any $\mathcal{F} \in \mathcal{TGL}(X)$ and let $\mathcal{O}$ be the $G^0$-orbit of $\mathcal{F}$.  Because $G^0$ is connected, this orbit either equals $\{ \mathcal{F} \}$ or is infinite.  We first deduce the conclusion of the main theorem in case $\mathcal{O}$ is infinite. 

\subsubsection{Warped product}
\label{section.warped.prod}

Consider the continuous map
\begin{eqnarray*}
\varphi & : & \mathcal{TGL}(X) \times X \rightarrow {\bf P}(TX)  \\
\varphi & : & ( \mathcal{F}, x ) \mapsto (x, (T \mathcal{F}_x)^{\perp})
\end{eqnarray*}

For each $x \in X$, the image $\varphi(\mathcal{O} \times \{ x \})$ is connected, so it is either infinite or just one point.  The set $D$ of all $x$ for which $| \varphi(\mathcal{O} \times \{ x \})| = 1$ is closed.  The complement $D^c \neq \emptyset$ because $\mathcal{O}$ is infinite.  For $x \in X$, let $C_x$ be the set of lightlike lines in $T_xX$ normal to leaves through $x$ of codimension-one, totally goedesic, lightlike hypersurfaces.  For all $x \in D^c$, the set $C_x$ is infinite.  

Now Theorem 1.1 of \cite{Zetgl2} applies to give an open set $U \subseteq D^c$ locally isometric to a warped product $N \times_h L$, where $N$ is Lorentzian of constant curvature, and $L$ is Riemannian.  For each $x \in U$, the subspace generated by $C_x$ equals $T_x N_x$, where $N_x$ is the $N$-fiber through $x$ (see the intermediate result \cite{Zetgl2} 3.3).  Since $X$ is the universal cover of a compact, real-analytic manifold, Theorem 1.2 of \cite{Zetgl2} implies that $X$ is a global warped product $N \times_h L$, and both $N$ and $L$ are complete.  

Because $G^0$ preserves the cone field $x \mapsto C_x$, it also preserves the $N$-foliation.  Then $G_1 = Isom^0(N) \lhd G^0$, so it is semisimple.  Since $X$ is contractible, $N$ and $L$ are, as well.  Then $N$ must be isometric to $\widetilde{AdS}^k$ for some $k$, and $G_1 \cong \widetilde{O}^0(2,k -1)$.  The assumption that $C_x \subset T_x N$ is infinite implies $k \geq 3$.

\subsubsection{Orbibundle}
\label{section.orbibundle}

Now it remains to show that $X \rightarrow G^0 \backslash X$ is a fiber bundle, and that $M$ is an orbibundle.  Let $G_2$ be the kernel of the homomorphism $G^0 \rightarrow Isom^0(N)$; it is semisimple, and $G^0 \cong G_1 \times G_2 \subseteq Isom(N) \times Isom(L)$.  The orbit $Y$ is isometric to $N \times L_2$ for a Riemannian submanifold $L_2$ of $L$, and $G_2$ is isomorphic to a subgroup of $Isom(L_2)$.  Clearly, $G_2(x)$ is compact for all $x \in X$.  We will show, using cohomological dimension, that $G_2(x) \cong K_2$ for all $x \in X$, where $K_2$ is a maximal compact subgroup of $G_2$.

Since $\Gamma_0$ acts properly discontinuously and cocompactly on $Y \cong N \times G_2/G_2(y)$, it is also properly discontinuous and cocompact on $N \times G_2/K_2$.  This latter space is contractible, so by Proposition \ref{basiccdq} (2),
$$\mbox{cd}_{\Q} \Gamma_0 = k + \mbox{dim}(G_2/K_2)$$

Next, the quotient $G_2 \backslash L$ can be identified with $G^0 \backslash X$.  Since $L$ is contractible and $G_2$ acts properly on it, either quotient is contractible by \cite{FW} 2.2.  We want to show that $\Gamma / \Gamma_0$ acts properly discontinuously on this quotient.  Suppose that a compact $\overline{C} \subset G_2 \backslash L$ is given.  The goal is to show that 
$$(\Gamma/\Gamma_0)_{\overline{C}} = \{ [ \gamma] \in \Gamma / \Gamma_0 : [\gamma] C \cap C \neq \emptyset \}$$

is finite.  

There is a compact $C \subset L_y \subset X$ projecting onto $\overline{C}$ by \ref{stratfacts} (1).  
Let $\mathcal{L}X$ be the bundle of Lorentz frames on $X$.  We may assume $C$ is small enough that $\mathcal{L}X \cong C \times O(1,n-1)$.  Let $A$ be the image of a continuous section of $\mathcal{L} X$ split along the product $X = N \times L$---that is, each frame in $A$ has the first $k$ vectors tangent to the $N$-foliation, and the succeeding vectors tangent to the $L$-foliation.  Let $B$ be the saturation $A \cdot (\Z_2 \times \Z_2 \times O(m))$, where $m = \mbox{dim} L$, and $\Z_2 \times \Z_2\subseteq O(1,k-1)$ acts transitively on orientation and time orientation of Lorentz frames along $N$; now $B$ is still compact.  Since $G$ acts properly on $\mathcal{L} X$ (see \cite{Gromov} 1.5.B or \cite{Ko} 3.2), the set $G_{A,B}$ is compact in $G$.  Because $N$ has constant curvature, $G_1 \cong Isom^0(N)$ is transitive on Lorentz frames along $N$, up to orienation and time orientation.  Then it is not hard to see 
$$G_{C,G^0 C} = G^0 \cdot G_{A,B} = G_{A,B} \cdot G^0$$

Then $G_{C,G^0 C}$ consists of finitely many components of $G$.  Now
$$(\Gamma/\Gamma_0)_{\overline{C}} = (\Gamma_{C, G^0C} \cdot \Gamma_0)/ \Gamma_0$$

Distinct $\Gamma_0$-cosets in $\Gamma$ occupy distinct components of $G$.  Then $\Gamma_{C,G^0 C}$ consists of finitely many cosets of $\Gamma_0$, and $(\Gamma/\Gamma_0)_{\overline{C}}$ is finite, as desired.

Now, as in Step 2 of Section \ref{section.proper.case}, 
$$\mbox{cd}_{\Q}(\Gamma/ \Gamma_0) \leq \mbox{dim}(G^0 \backslash X) = \mbox{dim}(G_2 \backslash L)$$

The inequality \ref{basiccdq} (3) gives
$$k + \mbox{dim}L \leq k + \mbox{dim} (G_2/K_2) + \mbox{dim}(G_2 \backslash L)$$

so $\mbox{dim} G_2(x) = \mbox{dim} K_2$, and $G_2(x)$ is conjugate in $G_2$ to $K_2$ for all $x$.  Then the quotient $\widetilde{Q} = G^0 \backslash X$ is a contractible manifold by \ref{stratfacts} (2). Since $\Gamma / \Gamma_0$ acts properly discontinuously here, $M$ is an orbibundle
$$ \Gamma_0 \backslash G^0 /H_0 \rightarrow M \rightarrow Q$$

The homogeneous space $G^0/H_0 \cong \widetilde{AdS}^k \times G_2/K_2$.

\subsubsection{Splitting}

From section \ref{section.warped.prod}, we have 
$$ X \cong \widetilde{AdS}^k \times_h L$$

where the warping function $h$ on $L$ is $G_2$-invariant.   The function $h$ descends to a function $h_1$ on $\widetilde{Q}$.  From the previous section, all $G_2$-orbits in $L$ are equivariantly diffeomorphic to $G_2 / K_2$.  As in the proper case, if $Z(G_2)$ is finite, we can define
\begin{eqnarray*}
 & \rho : X \rightarrow G_2 / K_2 & \\
 \rho(x) = [g] & \qquad \mbox{where} \qquad & gK_2 g^{-1} = G_2(x)
\end{eqnarray*}

This map factors through the projection to $L$.  As in the proper case, we can show that $L \cong G_2 / K_2 \times_{h_2} \widetilde{Q}$ for some $h_2 : \widetilde{Q} \rightarrow \mathcal{M}$, the moduli space of $G_2$-invariant Riemannian metrics on $G_2 / K_2$.   Now $h = (h_1, h_2)$ can be viewed as a function from $\widetilde{Q}$ to the moduli space of $G^0$-invariant Lorentz metrics on $ \widetilde{AdS}^k \times G_2 / K_2$.

\subsection{Nonproper case: no fixed point in $\mathcal{TGL}(X)$}
\label{tglx.big}

Now suppose, as above, that $G^0(y)$ is noncompact, so $\mathcal{TGL}(X) \neq \emptyset$, but every $G^0$-orbit in $\mathcal{TGL}(X)$ is a fixed point.  Then $G^0$ preserves a tgl foliation on $X$, so it preserves a lightlike line field on $X$.  We will show that this is impossible. 

First, we may assume that this lightlike line field along $Y$ is tangent to $Y$.  Suppose that $Y$ is either a fixed point or Riemannian.  Then the kernel of the restriction of $G^0$ to $Y$ contains a noncompact semisimple local factor $G_1$.  Recall that $C_y$ is the set of lightlike lines in $T_yX$ normal to codimension-one, totally geodesic, lightlike hypersurfaces through $y$.  Now $G_1$ acts on $C_y$ via the isotropy representation, and by assumption, it preserves an isotropic line in $C_y$, but this is impossible if $G_1$ is semisimple and noncompact.  Therefore, the orbit $Y$ is either Lorentzian or \emph{degenerate}---$T_yY^{\perp} \cap T_yY \neq {\bf 0}$ for all $y \in Y$.  If $Y$ is degenerate, then $G^0$ preserves the lightlike line field $T_yY^{\perp}$ along $Y$.  Suppose $Y$ is Lorentzian.  Now $G^0$ preserves the projections of the isotropic line field $y \mapsto C_y$ onto $TY$ and $(T Y)^{\perp}$.  If the second projection is nonzero, then the first is necessarily timelike.  But if $G^0(y)$ preserves a timelike vector in $T_yY$, then $G^0(y)$ must be compact, a contradiction.  Therefore, we may assume $G^0$ preserves a lightlike line field tangent to $Y$.

We first collect some facts about the isotropy representation.  We will show that it is either reductive or unimodular, in each case with a rather specific form.  Then we give standard examples of homogeneous spaces with each of these isotropy representations---the de Sitter space $dS^2$ and the light cone in Minkowski space---and show that they have no compact quotients.  Finally, we show that, in both the reductive and unimodular cases, the orbit $Y$ must essentially be one of these homogeneous spaces, contradicting that $\Gamma_0 \backslash Y$ is compact. 

\subsubsection{Properties of the isotropy respresentation}

Fix an isometric isomorphism of $T_y X$ with $\R^{1,n-1}$, determining an isomorphism $O(T_yX) \cong O(1,n-1)$.  Let $V$ be the image of $T_yY$ under this isomorphism, and let $k = \mbox{dim} V$.  Let $\Phi : G^0(y) \rightarrow O(1,n-1)$ be the resulting isotropy representation.  There is a filtration on $V$ preserved by $\Phi$.  The notation $U {\displaystyle \subset_i} V$ means $U$ is a subspace of $V$ with $\mbox{dim} (V/U) = i$.  The invariant filtration is
$$ {\bf 0} {\displaystyle \subset_1} V_0 {\displaystyle \subset_{k-1-i}} V_1 {\displaystyle \subset_i} V$$

where $i=0$ or $1$ depending on whether $V$ is degenerate or Lorentz.  The subspaces $V_0$ and $V_1$ are degenerate.  The quotient representation $V_1/V_0$ is orthogonal.  Because $\Phi$ preserves the isotropic line $V_0$ it descends to a quotient representation on $V_1/V_0$, which is orthogonal.  The image of $\Phi$ is conjugate in $O(1,n-1)$ to the minimal parabolic 
$$ P =  (M \times A) \ltimes U$$

where $U \cong \R^{n-2}$ is unipotent, $A \cong \R^*$, and $M \cong O(n-2)$, with the conjugation action of $M \times A$ on $U$ equivalent to the standard conformal representation of $O(n-2) \times \R^*$ on $\R^{n-2}$.  Denote by $\mathfrak{p}$ the Lie algebra of $P$, and by $\mathfrak{m}$, $\mathfrak{a}$, $\mathfrak{u}$, the subalgebras corresponding to $M$, $A,$ and $U$.   

Because $G^0$ acts properly and freely on the bundle of Lorentz frames of $X$, the isotropy representation is an injective, proper map.  By Corollary \ref{algstabs}, the image $\Phi(G^0(y))$ is algebraic.  Therefore, it decomposes 
$$ im(\Phi) \cong R^{\prime} \ltimes U^{\prime}$$

where $R^{\prime}$ is reductive and $U^{\prime}$ is unipotent (\cite{MostowSpl}).  Any unipotent subgroup of $P$ lies in $U$, so $U^{\prime} \subset U$.  The reductive complement $R^{\prime}$ is contained in a maximal reductive subgroup, which is then conjugate into $A \times M$.  Let $\mathfrak{r}^{\prime} \ltimes \mathfrak{u}^{\prime}$ be the corresponding Lie algebra decomposition.

Note that $T_yY$ can be identified with $\mathfrak{g}/ \mathfrak{g}(y)$ by the map $f_y$ as in Section \ref{kowalsky.argument}, and there is the relation
$$g_{*y} \circ f_y (B) =  f_y \circ \mbox{Ad}(g)(B)$$

for $B \in \mathfrak{g}$ and $g \in G^0(y)$.  In other words, $\Phi$ restricted to $V$ is equivalent to the representation $\overline{\mbox{Ad}}$ of $G^0(y)$ on $\mathfrak{g} / \mathfrak{g}(y)$ arising from the adjoint representation.  Let  $\varphi : \mathfrak{g}(y) \rightarrow \mathfrak{o}(1,n-1)$ be the Lie algebra representation tangent to $\Phi$ and $\ad$ be the representation tangent to $\overline{\mbox{Ad}}$. 

\begin{propn}
\label{filtration}
There is a filtration of $\mathfrak{g}$ invariant by the adjoint of $\mathfrak{g}(y)$:
$${\bf 0} \subset \mathfrak{g}(y) \subset_1 \mathfrak{s}(y) \subset_{k-1-i} \mathfrak{t}(y) \subset_i \mathfrak{g}$$

where $i=0$ or $1$ depending on whether $Y$ is degenerate or Lorentz.  The subspace $\mathfrak{s}(y)$ is a subalgebra.  The quotient representation for $\ad$ on $\mathfrak{t}(y) / \mathfrak{s}(y)$ is skew-symmetric.
\end{propn}

\begin{Pf}
The $\varphi$-invariant filtration ${\bf 0 } \subset V_0 \subset V_1 \subset V$ of $V$ corresponds to an $\ad$-invariant filtration of $\mathfrak{g} / \mathfrak{g}(y)$.  Lifting to $\mathfrak{g}$ gives the desired $\mbox{ad} (\mathfrak{g}(y))$-invariant filtration.  That $\mathfrak{s}(y)$ is a subalgebra follows from the facts that $[\mathfrak{g}(y),\mathfrak{s}(y)] \subset \mathfrak{s}(y)$ and $\mbox{dim}(\mathfrak{s}(y) / \mathfrak{g}(y)) = 1$.  Orthogonality of $\Phi$ on $V_1 / V_0$ implies $\varphi$ is skew-symmetric on $V_1 / V_0$; skew-symmetry of $\ad$ on $\mathfrak{t}(y) / \mathfrak{s}(y)$ follows.
\end{Pf}

Now we show that the image of $\Phi$ is either contained in $A \times M$ or $M \ltimes U$.

\begin{lmm}
The image of $\varphi$ is either reductive or consists of endomorphisms with no nonzero real eigenvalues.
\end{lmm}

\begin{Pf}
Suppose there is $B \in \mathfrak{g}(y)$ such that $\varphi(B)$ has nonzero eigenvalue $\lambda$ for some eigenvector ${\bf v} \in T_yX$.  The vector ${\bf v}$ is necessarily isotropic, and we may assume that ${\bf v} \in V_0$.  Otherwise, for any nonzero ${\bf w} \in V_0$, the inner product $<{\bf v}, {\bf w}> \neq 0$, which implies that $\varphi(B)$ has nonzero real eigenvalue on ${\bf w}$, as well.

Assume $\lambda > 0$; the case $\lambda < 0$ is similar.  We may assume $B \in \mathfrak{r}^{\prime}$.  The trace of $\left. \varphi(B) \right|_V$ is nonnegative and equals $0$ if and only if $V$ is Lorentz.  Correspondingly, the trace of $\ad(B)$ on $\mathfrak{g} / \mathfrak{g}(y)$ is nonnegative.

If $\varphi(B) \in \mathfrak{p}$ has eigenvalue $\lambda > 0$, then the adjoint $\mbox{ad}(\varphi(B))$ has no negative eigenvalues on $\mathfrak{p}$.  To simplify the argument, we will use that $\varphi(B) = B_1 + B_2$, where ${\bf 0} \neq B_1 \in \mathfrak{a}$ and $B_2 \in \mathfrak{m}$.  It is easy to see that $\mbox{ad}(B_1)$ has only real nonnegative eigenvalues on $\mathfrak{p}$.  All eigenvalues of $\mbox{ad}(B_2)$ are purely imaginary.  Since $\mbox{ad}(B_1)$ and $\mbox{ad}(B_2)$ are simultaneously diagonalizable, their sum $\mbox{ad}(\varphi(B))$ cannot have a negative eigenvalue.  

Now suppose that $im(\varphi)$ is not reductive, so $\mathfrak{u}^{\prime} \neq {\bf 0}$.  Let $m = dim (\mathfrak{u}^{\prime})$.  It is easy to compute that the trace of $\mbox{ad}(\varphi(B))$ on $\mathfrak{u}^{\prime}$ is $m \lambda$.  Since $\mbox{ad}(\varphi(B))$ has no negative eigenvalues, the trace of $\mbox{ad}(\varphi(B))$ on $im(\varphi) \subseteq \mathfrak{p}$ is positive.  Then the trace of $\mbox{ad}(B)$ on $\mathfrak{g}(y)$ is positive.

Finally, the trace of $\mbox{ad}(B)$ on $\mathfrak{g}$ is positive, which is impossible because $\mathfrak{g}$ is unimodular.
\end{Pf}

Now we have that $im(\Phi)$ is either a reductive subgroup of $A \times M$ or has the form $M^{\prime} \ltimes U^{\prime}$, where $M^{\prime} \subset M$ and $U^{\prime} \subset U$.

\subsubsection{Two examples with no compact quotient}
\label{examples.no.quotient}

{\bf Two-dimensional de Sitter space.  }  
\nopagebreak[3]

The $2$-dimensional de Sitter space $dS^2$ has isometry group $O(1,2)$ and isotropy $O(1,1)$, which has an index-two subgroup isomorphic to $\R^*$.  It is a well-known result of Calabi and Markus that no infinite subgroup of $O(1,2)$ acts properly on $dS^2$, so it has no compact quotient (\cite{CM}).
More generally, if $Y = dS^2 \times L$ for some Riemannian manifold $L$, then no subgroup of the product $O(1,2) \times Isom(L)$ acts properly discontinuously and cocompactly on $Y$; this is proved in \cite{Zetgl1} \S 15.1.

We will need an analogous result that also applies to the universal cover $\widetilde{dS}^2$.

\begin{propn}
\label{CMunivcover}
Let $G \cong \widetilde{O}^0(1,2)$, and $H$ be a connected Lie group.  There is no subgroup $\Gamma \subset G \times H$ acting properly discontinuously and cocompactly on $\widetilde{dS}^2 \times H$.
\end{propn}

\begin{Pf}
Let $K = Ad^{-1}(SO(2))$, where $SO(2)$ is a maximal compact subgroup of $Ad(G) \cong O^0(1,2)$.  Let $Z \cong \Z$ be the torsion-free factor of the center $Z(G)$.  Let $\overline{K}$ be a compact fundamental domain in $K$ for the $Z$-action with $\overline{K} = \overline{K}^{-1}$; for example, identifying $SO(2)$ with $S^1$ and $K$ with $\R$, we can take $\overline{K} = [-1/2,1/2]$.  Let $A$ be a maximal $\R$-split torus in $G$.  We have $G = K A K = Z \overline{G}$, where $ \overline{G} = \overline{K} A \overline{K}$.  For any $g \in \overline{G}$, 
$$ g \overline{K} \cap \overline{K} A \neq \emptyset$$

In other words, any $g$ in $\overline{G}$ takes the image $[ \overline{K} ]$ of $\overline{K}$, which is compact, in $\widetilde{dS}_2 = G/A$ to meet itself.

There is an isomorphism $G \cong \widetilde{SL}_2(\R)$, so $G$ acts on the real line, with $Z$ acting by integral translations.  
The $\emph{translation number}$
\begin{eqnarray*}
\tau & : & G \rightarrow \Z \cong Z \\
\tau & : & g \mapsto \lim \frac{g^n(0)}{n}
\end{eqnarray*}

is a continuous \emph{quasi-morphism} (see \cite{Ghys}): there exists $D > 0$ such that
$$| \tau(g g') - \tau(g) - \tau(g')| < D \qquad \mbox{for all} \ g,g' \in G$$
 
Note $\tau(\overline{K}) = [ -1/2, 1/2]$, and $\tau(A) = 0$.  Therefore if $g \in \overline{G}$, then $| \tau(g)| \leq 2D + 1$.  Also note that for $n \in \Z \cong Z$ and $g \in G$, then $\tau(ng) = n + \tau(g)$.

Now suppose that $\overline{C} \subset \widetilde{dS}^2 \times H$ is a compact fundamental domain for $\Gamma$.
  Denote by $\rho_1$ and $\rho_2$ the projections onto $\widetilde{dS}^2$ and $H$, respectively.  We may assume that the identity of $H$ is in $\rho_2(\overline{C}) = U$.  For $n \in Z$, let 
$$S_n = \{  (g, h) \in \Gamma \ : \ g \in n \overline{G}, \ h U \cap U \neq \emptyset \}$$

For a subset $L \subseteq G$, denote by $[L]$ its image in $\widetilde{dS}_2$.  Note that, for any $\gamma \in S_n$, the intersection 
$$\gamma ([ \overline{K}] \times U) \cap ( [n \overline{K}] \times U ) \neq \emptyset$$

Therefore, if $\Gamma$ acts properly discontinuously, then $|S_n| < \infty$ for each $n \in Z$.

On the other hand, we have $[ \overline{G}] \times U \subset \Gamma \cdot \overline{C}$.  Let $C$ be a compact lift of $\rho_1(\overline{C})$ to $G$.  Then we have  
$$\overline{G} \times U \subset \Gamma \cdot (C A \times U)$$  

The restriction of $| \tau |$ to $\overline{G} C A$ is bounded, so, for $|n|$ sufficiently large,
$$ n \overline{G} C A \cap \overline{G} = \emptyset$$

It follows that $\overline{G} \times U$ is contained in the union of finitely many $S_n \cdot (CA \times U)$, which is a union of finitely many translates $\gamma \cdot (C A \times U)$, which is impossible, because the image $[ \overline{G}] \times U$ is not compact.
\end{Pf}

{\bf The Minkowski light cone.  }
\nopagebreak[3]

A component of the light cone minus the origin in Minkowski space $\R^{1,k-1}$ is a degenerate orbit of $O^0(1,k-1)$, which we will momentarily denote by $G^0$.  The stabilizer of an isotropic vector is isomorphic to $M \ltimes U$, where $M, U \subset P$ are as above.  We will show that no subgroup of $G^0$ acts properly discontinuously and cocompactly on this orbit.

Suppose that $y$ is a point in the light cone and $\Gamma \subset G^0$ is a discrete subgroup such that $\Gamma \backslash G^0 / G^0(y)$ is a compact manifold.  Then $\Gamma \backslash G^0 / U$ is also compact; we may assume it is orientable.
Because $U$ is unimodular, the homogeneous space $G^0 / U$ has a $G^0$-invariant volume form (see \cite{Rag} I.1.4).  This form descends to $\Gamma \backslash G^0 / U$, where it has finite total volume.
The subgroup $A \cong \R^*$ of $P$ normalizes $U$, with generator $a$ acting by $Ad(a)(Y) = e^2Y$ for all $Y \in \mathfrak{u}$.  
Then $a$ acts on $\Gamma \backslash G^0 / U$ and scales the volume form by $1/e^{2(k-2)}$ at every point, which is impossible for a diffeomorphism of a compact manifold.

In the next section, we will show that, if $im(\Phi)$ is reductive, then $\widetilde{Y}$ is related, by proper $\widetilde{G}^0$-equivariant maps, to $\widetilde{dS}^2 \times H$, where $H$ is a connected Lie group.  In case $im(\Phi)$ is unimodular, we will show that there is a proper $G^0$-equivariant map $( O(1,k -1) / U) \times G_2 \rightarrow Y$, where $U$ is the unipotent radical of the minimal parabolic of $O(1,k -1)$, and $G_2$ is a local factor of $G^0$.  In both cases, no subgroup of $G^0$ can act properly discontinuously and cocompactly on $Y$.  Both cases involve studying the representation $\Phi$ and applying dynamical results from Section \ref{kowalsky.argument}.  

An element $B$ of $\mathfrak{g}$ is called \emph{nilpotent} if ad($B$) is nilpotent.  An element $B$ is \emph{semisimple} if ad$(B)$ is diagonalizable over $\C$, and $B$ is \emph{$\R$-split} if ad$(B)$ is diagonalizable over $\R$.

\subsubsection{Reductive case}

In this case, $im(\Phi) \subset A \times M$.  Because $G^0(y)$ is noncompact and $\Phi$ is proper, $im(\Phi)$ is not contained in $M$.  The image is fully reducible on $T_yX$; it decomposes as a product $A^{\prime} \times M^{\prime}$, where $M^{\prime}$ is compact, $A^{\prime}$ is one-dimensional, and $A^{\prime}$ has nontrivial character on $V_0$.  Let $\widehat{A} \times \widehat{M}$ be the corresponding decomposition of $G^0(y)$.  Properness of $\Phi$ implies $\widehat{M}$ is compact.  Continuity implies $\mbox{Ad}(a^n) \rightarrow \infty$ for all nontrivial $a \in \widehat{A}^0$: if $\mbox{Ad}(a^n)$ were bounded, then $\overline{\mbox{Ad}}(a^n)$ would be bounded, so $\Phi(a^n)$ would be bounded on $V$, a contradiction. 

\begin{propn}
\label{redfacts}
\  
\begin{enumerate}
\item{The restriction of the metric to $Y$ is Lorentzian, so $V_1 \neq V$.} 
\item{There is an $\overline{ad}$-invariant decomposition
$$\overline{\mathfrak{s}}_0(y) \oplus \overline{\mathfrak{s}}_1(y) \oplus \overline{\mathfrak{s}}_2(y)$$

of $\mathfrak{g}/\mathfrak{g}(y)$ corresponding to the filtration in Proposition \ref{filtration}.  
}

\item{The stabilizer subalgebra $\mathfrak{g}(y)$ contains no elements nilpotent in $\mathfrak{g}$; in particular, there are no root vectors of $\mathfrak{g}$ in $\mathfrak{g}(y)$.}
\end{enumerate}
\end{propn}

\begin{Pf}
\begin{enumerate}
\item{Let $B \in \widehat{\mathfrak{a}}$, and let $\lambda$ be the nonzero eigenvalue of $\varphi(B)$ on $V_0$, which we assume is positive.  If $V$ is degenerate, then the trace of $\varphi(B)$ on $V_1 = V$ is positive, so the trace of $\overline{\mbox{ad}}(B)$ on $\mathfrak{g} / \mathfrak{g}(y)$ is positive.  Now $B \in \mathfrak{z}(\mathfrak{g}(y))$, so the trace of $\mbox{ad}(B)$ on $\mathfrak{g}(y)$ is $0$.  Then the trace of $\mbox{ad}(B)$ on $\mathfrak{g}$ is positive, contradicting unimodularity of $\mathfrak{g}$.}

\item{Let $\mathfrak{s}(y) \subset \mathfrak{t}(y) \subset \mathfrak{g}$ be the $\mathfrak{g}(y)$-invariant subpaces in Proposition \ref{filtration}.  Let $\overline{\mathfrak{s}}_0(y)$ be the projection of $\mathfrak{s}(y)$ to $\overline{\mathfrak{g}} = \mathfrak{g}/ \mathfrak{g}(y)$.  Let $\overline{\mathfrak{t}}(y)$ be the projection of $\mathfrak{t}(y)$.  Because $\overline{\mbox{ad}}(\mathfrak{g}(y))$ is fully reducible, there is an invariant complement $\overline{\mathfrak{s}}_1(y)$ to $\overline{\mathfrak{s}}_0(y)$ in $\overline{\mathfrak{t}}(y)$.  Let $\overline{\mathfrak{s}}_2(y)$ be an invariant complement to $\overline{\mathfrak{t}}(y)$ in $\overline{\mathfrak{g}}$.}

\item{Suppose that $X \in \mathfrak{g}(y)$ is nilpotent.  Then $\overline{\mbox{ad}}(X)$ is nilpotent, so $\varphi(X)$ restricted to $V$ is nilpotent.  Because $im(\varphi)$ contains no nilpotent elements, $\varphi(X)$ is trivial on $V$.  By (1), the inner product on $V^{\perp} \subset R^{1,n-1}$ is positive definite, so $\varphi(X)$ is skew-symmetric and generates a precompact subgroup of $O(1,n-1)$.  Because $\Phi$ is proper, $X$ should generate a precompact subgroup of $G^0$, a contradiction unless $X = {\bf 0}$.}
\end{enumerate}
\end{Pf}

Now let $b \in \widehat{A}^0$, so $\mbox{Ad}(b^n) \rightarrow \infty$.  By Proposition \ref{kowarg}, there exists an $\R$-split element $B$ of $\mathfrak{g}$ and a root system such that 
$$\bigoplus_{\alpha(B) > 0} \mathfrak{g}_{\alpha} $$

is isotropic for the pullback inner product $<,>_y$ on $\mathfrak{g}$.
By Proposition \ref{redfacts} (3), this sum of root spaces does not meet $\mathfrak{g}(y)$.  Therefore, 
$$ \mbox{dim} \left( \bigoplus_{\alpha(B) > 0} \mathfrak{g}_{\alpha} \right) = 1$$

and some factor, say $\mathfrak{g}_1$, of $\mathfrak{g}$ is isomorphic to $\mathfrak{sl}_2(\R)$.    

Denote by $\mathcal{L}$ the null cone in $\mathfrak{g} / \mathfrak{g}(y) \cong V$.  For $b \in \widehat{A}^0$, the sequence $\overline{Ad}(b^n)$ has unique attracting and repelling fixed points, $p^+$ and $p^-$, respectively, in the projectivization ${\bf P}(\mathcal{L})$; these correspond to the nontrivial eigenvectors of $\overline{Ad}(b)$.  For $i = 1, \ldots, l$, denote by $\overline{\mathfrak{g}}_i$ the image of $\mathfrak{g}_i$ modulo $\mathfrak{g}(y)$; each such subspace is $\overline{Ad}(G(y))$-invariant.  Similarly for $X \in \mathfrak{g}$, denote by $\overline{X}$ its image in $\mathfrak{g} /\mathfrak{g}(y)$.  

Let $X_{\alpha}$ be a generator of the isotropic root space $\mathfrak{g}_{\alpha}$ above, so $\overline{X}_{\alpha} \in \mathcal{L} \cap \overline{\mathfrak{g}}_1$.  Then either the projectivization $[ \overline{X}_{\alpha}] = p^-$, or $[ \overline{Ad}(b^n)(\overline{X}_{\alpha})] \rightarrow p^+$.  In either case, one of $p^-, p^+$ is in $[\overline{\mathfrak{g}}_1]$, and $\overline{Ad}(b)$ has an eigenvector with nontrivial real eigenvalue in $\overline{\mathfrak{g}}_1$.

The fact that $\mathfrak{g}(y)$ contains no nilpotent elements (\ref{redfacts} (3)) implies that the Killing form $\kappa_1$ of $\mathfrak{g}_1 \cong \mathfrak{sl}_2(\R)$ is definite when restricted to $\mathfrak{g}_1(y)$.  The orthogonal is a $G^0(y)$-invariant complement $\mathfrak{l}_1$ that projects $Ad$-equivariantly and isomorphically onto $\overline{\mathfrak{g}}_1$.  The element of $\mathfrak{l}_1$ projecting to $\overline{X}_{\alpha}$ is an eigenvector for $Ad(b)$.  Because $\mathfrak{g}_1 \cong \mathfrak{sl}_2(\R)$, it follows that $\pi_1(b)$ is $\R$-split.  The other nontrivial eigenvector for $Ad(b)$ in $\mathfrak{l}_1$ projects to a different isotropic vector in $\mathfrak{g}/\mathfrak{g}(y)$ that is an eigenvector for $\overline{Ad}(b)$.  Then both $p^+$ and $p^-$ belong to $[ \overline{\mathfrak{g}}_1]$, and they are the images of nilpotent elements of $\mathfrak{g}_1$.

Now let $Y$ be any nilpotent in $\mathfrak{g}_i$ for $i > 1$.  By \ref{redfacts} (3), $\overline{Y} \neq {\bf 0}$.  If $Y \notin \overline{s}_1(y)$, then $[\overline{Ad}(b^n)(\overline{Y})]$ converges in ${\bf P}(V)$ to $p^+$.  Then $p^+$ would be the image of a nilpotent element from $\mathfrak{g}_1$ and another from $\mathfrak{g}_i$.  In the span of these two would be a nilpotent element of $\mathfrak{g}(y)$, a contradiction.  Therefore, $\overline{\mathfrak{g}}_i \subseteq \overline{\mathfrak{s}}_1(y)$, and $\overline{Ad}(b^n)$ is bounded on $\overline{\mathfrak{g}}_i$ for all $i > 1$.

Since $\pi_1(b) \neq 1$, the intersection $G_i \cap \widehat{A}^0 = 1$ for all $i > 1$.  It follows that $\mathfrak{g}_i(y) \subseteq \mathfrak{m}'$, so it is definite for the restriction of the Killing form $\kappa_i$ of $\mathfrak{g}_i$, and $Ad(G^0(y))$ is bounded on $\mathfrak{g}_i(y)$.  As above, the orthogonal of $\mathfrak{g}_i(y)$ is an $Ad(G^0(y))$-invariant complement in $\mathfrak{g}_i$, which projects equivariantly and isomorphically to $\overline{\mathfrak{g}}_i$.  Now for all $i > 1$, the adjoint $Ad(G^0(y))$ is bounded on $\mathfrak{g}_i$, which implies that $Ad(\pi_i(G^0(y)))$ is precompact.  

Because $G^0$ preserves a Lorentz metric on $Y \cong G^0/G^0(y)$, any element of $Z(G^0) \cap G^0(y)$ would have trivial derivative along $Y$ at $y$, and would therefore be trivial on $Y$.  For the contradiction we are about to obtain, we may assume that $G^0$ acts faithfully on $Y$, so we may assume that $Z(G^0) \cap G^0(y) = 1$.  Therefore, $\pi_i(G^0(y))$ is precompact; let $K_i$ denote the compact closure.  Let $K = K_2 \times \cdots \times K_l$. 

We have already established that the projection $\pi_1(\widehat{A})$ contains a nontrivial $\R$-split element.  Because any other element of $\pi_1(G^0(y))$ must centralize this one, we conclude that $\pi_1(G^0(y))$ is isogenous to the maximal $\R$-split subgroup $A_1$ of $G_1$.

Therefore, $G^0(y) \subseteq A_1 \times K$, and there is a $G^0$-equivariant proper map
$$Y \cong G^0/G^0(y) \rightarrow G^0/(A_1 \times K)$$

so $\Gamma_0$ acts properly and cocompactly on both spaces.  Then $\Gamma_0$ acts properly and cocompactly on $ G/ A_1 \times H$, where
$H \cong G_2 \times \cdots \times G_l$.  Taking the preimage of $\Gamma_0$ in $\widetilde{G}^0$ if necessary, we obtain a subgroup of $\widetilde{O}^0(1,2) \times \widetilde{H}$ acting properly and cocompactly on $\widetilde{dS}^2 \times \widetilde{H}$, contradicting Proposition \ref{CMunivcover}.

\subsubsection{Unimodular case}

Now assume $im(\Phi) = M^{\prime} \ltimes U^{\prime}$ with $M^{\prime}$ compact and $U^{\prime}$ unipotent.  First we collect some algebraic facts for this case.

\begin{propn}
\label{unifacts}
Let $B \in \mathfrak{g}(y)$.
\begin{enumerate}
\item{$B$ is not $\R$-split.}

\item{If $\varphi(B)$ is nilpotent, then $B$ is nilpotent.}

\item{If $B$ is nilpotent, then on the filtration in Proposition \ref{filtration}, $\overline{\mbox{ad}}(B)$ carries each subspace to the next.  In other words, $\overline{\mbox{ad}}(B)$ is trivial on each factor of the associated graded space.}
 
\item{If $\varphi(B)$ is nilpotent and $\mathfrak{g} \neq \mathfrak{t}(y)$, then $\overline{\mbox{ad}}(B)$ has nilpotence order $3$.}   


\end{enumerate}
\end{propn}

\begin{Pf}
\begin{enumerate}
\item{Suppose $B$ is $\R$-split.  Let $\alpha$ be a root with $\alpha(B) \neq 0$ and $X_{\alpha}, X_{- \alpha}$ generators of the corresponding root spaces. The elements $X_{\alpha}, X_{- \alpha}$ generate a subalgebra of $\mathfrak{g}$ isomorphic to $\mathfrak{sl}_2(\R)$.  Because $\varphi(B)$ can have no eigenvectors with nonzero real eigenvalue, $X_{\alpha}$ and $X_{- \alpha}$ are both contained in $\mathfrak{g}(y)$.  Then $\mathfrak{g}(y) \subset \mathfrak{p}$ contains a subalgebra isomorphic to $\mathfrak{sl}_2(\R)$, a contradiction.  }

\item{If $\varphi(B)$ is nilpotent, then $\varphi(B) \in \mathfrak{u}^{\prime}$. Then $\overline{\mbox{ad}} (B)$ is nilpotent and $\mbox{ad}(B)$ is nilpotent on $\mathfrak{g}(y)$, which implies nilpotence of $B$.}

\item{If $B$ is nilpotent, then $\mbox{ad}(B)$ is trivial on both $\mathfrak{g} / \mathfrak{t}(y)$ and $\mathfrak{s}(y) / \mathfrak{g}(y)$, because they are both at most one-dimensional.  Because $\mbox{ad}(B)$ is skew-symmetric on $\mathfrak{t}(y) / \mathfrak{s}(y)$, this representation is also trivial.}

\item{If $\varphi(B)$ is nilpotent and $\mathfrak{g} \neq \mathfrak{t}(y)$, then $V$ is Lorentzian and the inner product on $V^{\perp}$ is positive definite, so $\varphi(B)$ is trivial on $V^{\perp}$.  By injectivity of $\varphi$, the restriction of $\varphi(B)$ to $V$ is nontrivial, so $\overline{\mbox{ad}}(B)$ is nontrivial.  By item (3), $\overline{\mbox{ad}}(B)$ has nilpotence order at most $3$.  Let $W \in \mathfrak{g} \backslash \mathfrak{t}(y)$.  We will show that $\mbox{ad}^2(B)(W) \notin \mathfrak{g}(y)$.  

Denote by $<,>$ the pullback of the inner product from $T_yX$ to $\mathfrak{g}$.  For any $W,Z \in \mathfrak{g}$
$$< \mbox{ad}(B)(W), Z> + <W, \mbox{ad}(B)(Z)> = 0$$

First we will show that $\mbox{ad}(B)(W) \notin \mathfrak{s}(y)$.  Suppose it is.  For any $Z \in \mathfrak{s}(y) \backslash \mathfrak{g}(y)$, the inner product $<W,Z> \neq 0$.  The identity
$$ < \mbox{ad}(B)(W), W> + <W , \mbox{ad}(B)(W) > = 2<\mbox{ad}(B)(W),W> = 0$$

implies $\mbox{ad}(B)(W) \in \mathfrak{g}(y)$.  Now $\mbox{ad}(B)(\mathfrak{t}(y)) \subseteq \mathfrak{g}(y)$ would imply $\overline{\mbox{ad}}(B)$ is trivial, which cannot be.  Then there must be some $Z \in \mathfrak{t}(y)$ such that $\mbox{ad}(B)(Z) \in \mathfrak{s}(y) \backslash \mathfrak{g}(y)$.  Then 
$$ < \mbox{ad}(B)(W), Z> = - < W, \mbox{ad}(B)(Z)> \neq 0$$

contradicting that $\mbox{ad}(B)(W) \in \mathfrak{s}(y)$.  

Now $\mbox{ad}(B)(W)$ must be in $\mathfrak{t}(y) \backslash \mathfrak{s}(y)$, so
$$< \mbox{ad}(B)^2 (W) ,W> = - < \mbox{ad}(B)(W), \mbox{ad}(B)(W) > \neq 0$$

which implies $\mbox{ad}(B)^2 (W) \notin \mathfrak{g}(y)$, as desired.
}

\end{enumerate}
\end{Pf}

 Let $\widehat{M} \ltimes \widehat{U}$ be the decomposition of $G^0(y)$ corresponding to $im(\Phi) = M^{\prime} \ltimes U^{\prime}$.  Again, because $\Phi$ is proper, $\widehat{M}$ is compact.  Let $\widehat{\mathfrak{m}}$ and $\widehat{\mathfrak{u}}$ be the corresponding subalgebras of $\mathfrak{g}$.  From item (2) above, $\widehat{u}$ consists of nilpotent elements.
Let $J$ be the set of $i$ such that $\pi_i(\widehat{\mathfrak{u}}) \neq {\bf 0}$. 

We can show by induction that there exists $X \in \widehat{\mathfrak{u}}$ such that $\pi_i(X) \neq {\bf 0}$ if and only if $i \in J$.  Let $i_1, \ldots, i_k$ be some order on the elements of $J$.  Clearly, there is some $X_1 \in \widehat{\mathfrak{u}}$ such that $\pi_{i_1}(X_1) \neq  {\bf 0}$.  Suppose $X_m \in \widehat{\mathfrak{u}}$ is such that $\pi_{i_j}(X_m) \neq {\bf 0}$ for all $j \leq m$.  There exists $Y_m \in \widehat{\mathfrak{u}}$ such that $\pi_{i_{m+1}}(Y_m) \neq {\bf 0}$.  For some real number $c$, the element $X_{m+1} = X_m + cY_m$ will have $\pi_{i_j}(X_{m+1}) \neq {\bf 0}$ for all $j \leq m+1$.  Write this element $X = \sum_{i \in J} X_i$ with $X_i \in \mathfrak{g}_i$.  Note that nilpotence of $X$ implies nilpotence of each $X_i$.

The Jacobson-Morozov theorem (see \cite{Helg} IX.7.4) yields, for each $i \in J$, an $\R$-split element $A_i \in \mathfrak{g}_i$ and a nilpotent element $Y_i \in \mathfrak{g}_i$ such that 
$$ [A_i,X_i] = 2X_i, \qquad [A_i,Y_i] = -2Y_i, \qquad \mbox{and} \ [X_i,Y_i] = A_i$$

Then the elements $A = \sum_i A_i$ and $Y = \sum_i Y_i$ satisfy
$$[A,X] = 2X, \qquad [A,Y] = -2Y, \qquad \mbox{and} \ [X,Y] = A$$

The subalgebra generated by $X,A,$ and $Y$ is isomorphic to $\mathfrak{sl}_2$.  Let $L$ be the corresponding subgroup of $G^0$.  
The adjoint of $\mathfrak{g}(y)$ is trivial on $\mathfrak{s}(y)/\mathfrak{g}(y)$ by \ref{unifacts} (3), so $G^0$ preserves a vector field tangent to the invariant isotropic line field along $Y$.  Now Proposition \ref{kowarg} and Remark \ref{kak.remark} for $g_n = e^{nX}$ give some $k \in L$ such that, for $A^{\prime} = \mbox{Ad}(k)(A)$, 
$$\bigoplus_{\alpha(A^{\prime}) > 0} \mathfrak{g}_{\alpha} \subset \mathfrak{s}(y)$$

Let $X^{\prime} = \mbox{Ad}(k)(X) \in \mathfrak{s}(y)$.  We may assume $k \in PSL_2(\R)$. For

$$
k = 
\left( 
\begin{array}{cc} 
           \cos \theta & - \sin \theta \\
           \sin \theta & \cos \theta
\end{array}
\right)
$$

$\mbox{and}$

$$
X = 
\left( 
\begin{array}{cc}
           0 & 1 \\
           0 & 0 
\end{array}
\right)
$$

the bracket
$$[X,X^{\prime}] = \left( \begin{array}{cc} 
                         - \sin ^2 \theta & 2 \cos \theta \sin \theta \\
                         0 &  \sin ^2 \theta
\end{array}
\right)
$$

This bracket belongs to $\mathfrak{g}(y)$, which contains no $\R$-split elements (\ref{unifacts} (3) and (1)).  Then $\sin \theta$ must be 0, so $\mbox{Ad}(k)$ is trivial, and
$$\bigoplus_{\alpha(A) > 0} \mathfrak{g}_{\alpha} \subset \mathfrak{s}(y)$$

Now we will show that this sum of root spaces is in fact contained in $\mathfrak{g}(y)$.  For $Y$ the negative root vector as above, $\mbox{ad}^2(X)(Y) \in \mathfrak{g}(y)$, so $\overline{\mbox{ad}}(X)$ has order less than $3$ on the corresponding element of $\mathfrak{g} / \mathfrak{g}(y)$.  Then $Y \in \mathfrak{t}(y)$ by Proposition \ref{unifacts} (4).  Then $[X,Y] = A \in \mathfrak{s}(y)$ by \ref{unifacts} (3), but $A$ cannot be in $\mathfrak{g}(y)$ by \ref{unifacts} (1), so
$$ \mathfrak{s}(y) = \R A + \mathfrak{g}(y)$$

Now suppose $\alpha(A) > 0$ and let $X^{\prime}$ be an arbitrary element of $\mathfrak{g}_{\alpha} \subseteq \mathfrak{s}(y)$.  Since $[\mathfrak{s}(y), \mathfrak{s}(y)] \subseteq \mathfrak{g}(y)$, the bracket $[A,X^{\prime}] = \alpha(A) X^{\prime} \in \mathfrak{g}(y)$. Therefore, $\oplus_{\alpha(A) > 0} \mathfrak{g}_{\alpha} \subseteq \mathfrak{g}(y)$, as desired; in particular, $X_i \in \mathfrak{g}(y)$ for all $i \in J$.

Next we will show that $|J|=1$.  As above, Proposition \ref{unifacts} implies $Y_i \in \mathfrak{t}(y)$ and $A_i \in \mathfrak{s}(y)$ for all $i \in J$.  If $|J| > 1$, then, for one $i \in J$ and some nonzero $c \in \R$, the difference $cA - A_i$ is a nontrivial $\R$-split element of $\mathfrak{g}(y)$, contradicting \ref{unifacts} (1).  

Now $G^0(y)$ has precompact projection on all local factors but one, say $G_1$.  By Kowalsky's Theorem (\cite{Kowthesis}), $\mathfrak{g}_1 \cong \mathfrak{o}(2,k)$, for some $k \geq 3$, or $\mathfrak{o}(1,k)$, for some $k \geq 2$.  We will deduce that $\mathfrak{g}_1$ must be the latter, and that $G_1(y)$ is as in the Minkowski light cone.  The subspaces $\mathfrak{s}_1(y)$ and $\mathfrak{t}_1(y)$ will denote the intersections $\mathfrak{g}_1 \cap \mathfrak{s}(y)$ and $\mathfrak{g}_1 \cap \mathfrak{t}(y)$, respectively, below.

{\bf Step 1: } $\mathfrak{g}_1 \cong \mathfrak{o}(1,2)$ implies $G_1y$ degenerate.
\nopagebreak[3]

If $\mathfrak{g}_1 \cong \mathfrak{o}(1,2)$, then it is generated by $X,A$, and $Y$ from above.  Recall $X \in \mathfrak{g}_1(y); A \in \mathfrak{s}_1(y)$; and $Y \in \mathfrak{t}_1(y)$.  Then $\mathfrak{g}_1/\mathfrak{g}_1(y)$ is $2$-dimensional and degenerate with respect to the inner product pulled back from $T_yX$; therefore, the orbit $G_1y$ is also degenerate.

{\bf Step 2: } $G_1 y$ is degenerate in general.
\nopagebreak[3]

Assume that $\mathfrak{g}_1$ is not isomorphic to $\mathfrak{o}(1,2)$, and suppose that the orbit $G_1 y \subseteq Y$ is of Lorentzian type.  Then Theorem 1.5 of \cite{ADZ} gives that $G_1 y$ is equivariantly homothetic, up to covers, to $dS^k$ or $AdS^k$ for some $k \geq 3$; in either case, $\mathfrak{g}_1(y)$ would be semisimple, a contradiction.


 
{\bf Step 3: } Case $\mathfrak{g}_1 \cong \mathfrak{o}(2,k)$.
\nopagebreak[3]

Now suppose $\mathfrak{g}_1 \cong \mathfrak{o}(2,k)$ for some $k \geq 3$.   Let $\Delta$ be a root system of $\mathfrak{g}$ as above.  Let $A \in \mathfrak{s}(y)$ be as above.  Let $\alpha \in \Delta$ be such that $\alpha(A) = 2$. Let $X \in \mathfrak{g}_{\alpha} \cap \mathfrak{g}(y)$.  The root system of $\mathfrak{o}(2,k)$ is generated by two simple roots, $\beta$ and $\gamma$.  The root spaces for $\beta$ and $\gamma$ are each $(k -2)$-dimensional.  The other positive roots are $\beta - \gamma$ and $\beta + \gamma$, with one-dimensional root spaces.  

First suppose $\alpha = \beta$, so $X \in \mathfrak{g}_{\beta} \subset \mathfrak{g}(y)$.  Let $L$ be a generator of $\mathfrak{g}_{- \beta - \gamma}$.  For any such $X$ and $L$, the adjoint $\mbox{ad}^2(X)(L) \neq {\bf 0}$.  Since the orbit $G_1y$ is degenerate, $L \in \mathfrak{t}(y)$, and $\mbox{ad}(X)(L) \in \mathfrak{s}(y)$.  Let $W = \mbox{ad}(X)(L) \in \mathfrak{g}_{- \gamma}$.  Any nilpotent subalgebra of $\mathfrak{g}_1(y)$ is abelian, so $W \in \mathfrak{s}(y) \backslash \mathfrak{g}(y)$.  Then $cW - A \in \mathfrak{g}_1(y)$ for some nonzero $c \in \R$.  But now $L \in \mathfrak{t}(y) \backslash \mathfrak{s}(y)$ would be an eigenvector for this element with nonzero real eigenvalue, contradicting that $\mathfrak{g}_1(y)(L) \subseteq \mathfrak{s}(y)$.


We conclude that $X$ cannot be in $\mathfrak{g}_{\beta}$.  The same argument shows $X$ cannot be in $\mathfrak{g}_{\gamma}$; in fact, $\mathfrak{g}_1(y) \cap \mathfrak{g}_{\omega}$ must be ${\bf 0}$ for $\omega =  \pm \beta, \pm \gamma$.  

Now suppose that $\alpha = \beta \pm \gamma$, so either $(\beta + \gamma)(A)$ or $(\beta - \gamma)(A)$ equals $2$.  Then one of $\beta(A)$ or $\gamma(A)$ is nonzero, which again implies that one of $\mathfrak{g}_{\pm \beta}, \mathfrak{g}_{\pm \gamma}$ is in $\mathfrak{g}(y)$, a contradiction.





{\bf $G_1 y$ is the Minkowski light cone}
\nopagebreak[3]

Now we have that $\mathfrak{g}_1 \cong \mathfrak{o}(1,k)$ for some $k \geq 2$.  Let $\alpha$ be the positive root of $\mathfrak{g}_1$ with $\alpha(A) = 2$.  From above, $\mathfrak{g}_{\alpha} \subset \widehat{\mathfrak{u}}$. Since this root space is a maximal abelian subalgebra of nilpotent elements in $\mathfrak{g}_1$, this containment is equality by \ref{unifacts} (2).

Now $e^A = a$ normalizes $\widehat{U}$.  There is a proper equivariant map
$$ G^0 / \widehat{U} \rightarrow G^0/(\widehat{M} \ltimes \widehat{U}) \cong Y$$

so no subgroup of $G^0$ acts properly discontinuously and cocompactly on $Y$, as in Section \ref{examples.no.quotient}.

\section{Appendix: Uniformly Lipschitz foliations}
\label{appendix.lipschitz.foliations}

The purpose of this section is to provide the proof of Proposition \ref{uniformly.lipschitz}: 

\emph{
Let $X$ be the universal cover of a compact manifold $M$.  Let $\nabla$ be a smooth connection and $\sigma$ a smooth Riemannian metric, both lifted from $M$.   For any $r>0$, there exist $C, \delta > 0$ such that any radius-$r$, codimension-one geodesic lamination $(X^{\prime},f)$ on $X$ is $(C,\delta)$-Lipschitz:  any $x,y \in X^{\prime}$ with $d_{\sigma}(x,y) < \delta$ are connected by a unique $\nabla$-geodesic $\gamma$, and 
$$ \angle_{\sigma} (P_{\gamma} f(x), f(y)) \leq C \cdot d_{\sigma}(x,y)$$
}

Refer to Definition \ref{defn.geod.lamn} for \emph{codimension-one geodesic lamination}.

\emph{Notation}:  Let $\mbox{dim} X = k$.  Since $X$ is orientable, the unit normal $N(H)$ of an oriented hyperplane $H \subset T_xX$ is well-defined.  For the remainder of this section, all metric notions, such as the distance $d$, norm $| \cdot |$, length $l$, unit normal $N$, and angle $ \angle$, always refer to $\sigma$ below.  In any metric space below, $B(x, \delta)$ and $D(x, \delta)$ denote the ball and the disk, respectively, of radius $\delta$ around $x$.  For a Riemannian manifold $X$, the ball bundle $B^{\delta}X$ is the union of the $\delta$-balls about ${\bf 0}$ in each tangent space; $D^{\delta}X$ is the analogous disk bundle.  The bundle of spheres of radius $\delta$ will be denoted $T^{\delta}X$.  Affine notions below, such as geodesics, parallel transport $P$, and the exponential map $exp$, refer to $\nabla$, unless indicated otherwise.  For a smooth curve $\gamma$, parallel transport along $\gamma$ for time $t$ will be denoted $P^t_{\gamma}$.  If $\gamma$ is a continuous, piecewise smooth curve with explicit compact domain, the $P_{\gamma}$ is parallel transport along $\gamma$ from the initial to terminal point.

For two subspaces $W$ and $W^{\prime}$ of a finite-dimensional vector space $V$ with positive-definite inner product, the angle between $W$ and $W^{\prime}$ will mean
$$ \angle(W, W^{\prime}) = d_H(W \cap S({\bf 0},1), W' \cap S({\bf 0}, 1))$$

where $S({\bf 0},1)$ is the unit sphere in $V$, and $d_H$ is the Hausdorff distance.
We first establish several lemmas that will be needed for the proof of the proposition.

\begin{defn}
A neighborhood $B$ of a point $x \in X$ is \emph{hypersurface foliated based at $x$} if, for any oriented hyperplane $H \subset T_xX$, there is a neighborhood $W_{H,x}$ of $({\bf 0},0)$ in $H \times \R$ such that
$$\varphi_{H,x}  :  ({\bf p}, t) \mapsto exp_{\gamma_H(t)}(P^t_{\gamma_H} {\bf p})$$

where $\gamma_H(t) = exp_x(t N(H))$, is a diffeomorphism from $W_{H,x}$ onto $B$.

An open set $B \subset X$ is \emph{hypersurface foliated} if it is hypersurface foliated based at each $x \in B$.
\end{defn}

\begin{lmm}
\label{hypersurface.foliated}
There exists $\delta > 0$ such that for every $x \in X$, the ball $B(x,\delta)$ is hypersurface foliated.
\end{lmm}

\begin{Pf}
Because $\sigma$ and $\nabla$ are lifted from $M$, it suffices to find $\delta$ such that each $B(x,\delta)$ in the compact quotient $M$ is hypersurface foliated, assuming we choose $\delta$ sufficiently small that each $B(x, \delta)$ in $X$ projects diffeomorphically to its image in $M$.  We may assume $M$ is oriented.

Let $\Omega = SO(k)$, and let $A \subset M$ be an open subset equipped with a trivialization $\tau : A \times \Omega \rightarrow \mathcal{O}A$, where $\mathcal{O}A$ is the bundle of positively-oriented orthonormal frames (with respect to $\sigma$) on $A$. Given $x \in A$ and $\omega \in \Omega$, write 
$$ \tau(x,\omega) = (\tau_1(x,\omega), \ldots, \tau_k(x,\omega))$$

for the elements of the basis of $T_xX$ given by $\tau$.

Given $x \in A$ and $\omega \in \Omega$, there is a neighborhood $Y_{\omega,x}$ of $({\bf 0}, 0)$ in $\R^{k-1} \times \R$ on which the map $$ \varphi_{\omega,x} :  ({\bf p}, t) \mapsto exp_{\gamma_{\omega}(t)}(P^t_{\gamma_{\omega}} (\sum_{i=1}^{k-1} p_i \tau_i(x,\omega))) \in M$$

where $\gamma_{\omega}(t) = exp_x(t \tau_k(x,\omega))$, is defined.  Note that for 
$$H = span \{ \tau_1(x,\omega), \ldots, \tau_{k-1}(x,\omega) \}$$

 the vector $\tau_k(x,\omega) = N(H)$, and $\varphi_{\omega,x} = \varphi_{H,x}$ as in the definition of hypersurface foliated, under the isomorphism of $\R^{k-1}$ with $H$ given by $\tau(x,\omega)$.  Now consider the smooth map, defined on 
$$Y = \cup_{\omega \in \Omega, x \in A} ( \{ \omega \} \times \{ x \} \times Y_{\omega, x} ) \subseteq \Omega \times A \times \R^{k-1} \times \R$$ 

given by
$$\Phi  : (\omega,x,{\bf p}, t) \mapsto (\omega, x,\varphi_{\omega,x}({\bf p},t)) \in \Omega \times A \times M$$

The restriction of $\Phi$ to any $\{ \omega \} \times \{ x \} \times Y_{\omega,x}$ agrees with $\varphi_{\omega,x}$.  For each $\omega$ and $x$, the differential $D_{({\bf 0}, 0)} \varphi_{\omega,x}$ is an isomorphism.  Therefore, at any $(\omega,x) \in \Omega \times A$, the differential $D_{(\omega,x, {\bf 0}, 0)} \Phi$ is an isomorphism.  The inverse function theorem yields a neighborhood $U$ of $(\omega, x, x)$ such that $\Phi$ is a diffeomorphism from $\Phi^{-1}(U)$ onto $U$.  There exist $\delta_{\omega,x} > 0$ and a neighborhood $V_{\omega,x}$ of $\omega$ such that $V_{\omega,x} \times B(x,\delta_{\omega,x}) \times B(x,\delta_{\omega,x}) \subset U$.  Therefore, for every $\omega^{\prime} \in V_{\omega,x}$ and $y \in B(x,\delta_{\omega,x})$, the map $\varphi_{\omega^{\prime},y}$ is a diffeomorphism from some neighborhood of $({\bf 0},0)$ in $\R^{k-1} \times \R$ onto $B(x, \delta_{\omega,x})$.

For each $x \in M$, the collection of all $V_{\omega,x}$ is a cover of $\Omega$.  Let $V_{\omega_1,x}, \ldots, V_{\omega_N,x}$ be a finite subcover.  Let $\delta_x = min \{ \delta_{\omega_1,x}, \ldots, \delta_{\omega_N,x} \}$.  The balls $B(x,\delta_x)$ are all hypersurface foliated and form a cover of $M$.  Let $\delta$ be the Lebesgue number of this cover.  Then for all $x \in M$, the ball $B(x, \delta)$ is hypersurface foliated.
\end{Pf}

\begin{defn}
A neighborhood $B$ of a point $x \in X$ is \emph{normally foliated based at $x$} if, for any oriented hyperplane $H \subset T_xX$, there is a neighborhood $W_{H,x}$ of $({\bf 0}, 0)$ in $H \times \R$ such that
$$\varphi_{H,x} : ({\bf p}, t) \mapsto exp_{\gamma_{\bf p}(1)}(t N(P^1_{\gamma_{\bf p}} H))$$

where $\gamma_{\bf p}(t) = exp_x(t {\bf p})$, is a diffeomorphism from $W_{H,x}$ onto $B$.
\end{defn}

\begin{lmm}
\label{normally.foliated}
There exists $\delta > 0$ such that, for any $x \in X$, the ball $B(x, \delta)$ is normally foliated based at $x$. Further, for any $\eta > 0$, it is possible to choose $\delta$ such that, for all $H,x$, the inverse image $\varphi_{H,x}^{-1}(B(x,\delta)) \subseteq B({\bf 0}, \eta)$.  
\end{lmm}

\begin{Pf}
Let $\Omega$ and $A$ be as in the proof of \ref{hypersurface.foliated}.  For $\omega \in \Omega$ and $x \in A$, there is a neighborhood $Y_{\omega,x}$ of $({\bf 0}, 0)$ in $\R^{k-1} \times \R$ on which 
$$\varphi_{\omega,x}  : ({\bf p}, t) \mapsto exp_{\gamma_{\bf p}(1)}(t N (P^1_{\gamma_{\bf p}} span \{ \tau_1(x,\omega), \ldots, \tau_{k-1}(\omega,x) \} )) \in M $$

where $\gamma_{\bf p}(t) = exp_x(t {\bf p})$, is defined.  We require that $Y_{\omega,x} \subseteq B({\bf 0}, \eta)$ for all $\omega, x$.

On the union
$$ Y = \cup_{\omega \in \Omega, x \in A} ( \{ \omega \} \times \{ x \} \times Y_{\omega,x} ) \subseteq \Omega \times A \times B({\bf 0}, \eta) \subset \Omega \times A \times \R^{k-1} \times \R$$

define
$$\Phi : (\omega,x,{\bf p}, t) \mapsto (\omega, x,\varphi_{\omega,x}({\bf p},t)) \in \Omega \times A \times M$$

Now following the proof of \ref{hypersurface.foliated}, we can actually find $\delta > 0$ such that for all $x \in X$, the ball $B(x, \delta)$ is normally foliated based at each $ y \in B(x, \delta)$; in particular, $B(x,\delta)$ is normally foliated based at $x$.  Since the domain $Y$ is contained in $\Omega \times A \times B({\bf 0}, \eta)$, the second claim of the lemma is satisfied.
\end{Pf}

The following three lemmas establish bounds or Lipschitz relations between the metric $\sigma$ and the connection $\nabla$.  Let $\epsilon$ be such that, for every $x \in M$, the ball $B(x,2 \epsilon)$ is a normal neighborhood of each of its points for both $\nabla$ and the Levi-Civita connection for $\sigma$.  Then each $B(x, 2 \epsilon)$, $x \in X$, has the same property.

\begin{lmm}
\label{one.nabla.vs.sigma}
There exists $C > 0$ such that, for any $x,y \in X$ with $d(x,y) < \epsilon$,
$$ l(\gamma) \leq C \cdot d(x,y)$$

 where $\gamma$ is the $\nabla$-geodesic connecting $x$ to $y$.
\end{lmm}

\begin{lmm}
\label{two.nabla.vs.sigma}
For any $0 < \theta \leq \pi$, there exists $0 < \delta \leq \epsilon$ such that, for all $x,y \in X$ with $d(x,y) < \delta$,
$$| \angle(\gamma_{\sigma}^{\prime}(0), \gamma_{\nabla}^{\prime}(0))| \leq \theta$$

where $\gamma_{\sigma}$ and $\gamma_{\nabla}$ are the $\sigma$- and $\nabla$-geodesics, respectively, between $x$ and $y$.
\end{lmm}

\begin{lmm}
\label{three.nabla.vs.sigma}
There exists $C > 0$ such that, for any $\nabla$-geodesic triangle $T = \Delta(xyz)$ with $max \{ d(x,y), d(y,z), d(z,x) \} < \epsilon$, and any ${\bf u} \in T_xX$, 
$$| \angle({\bf u},P_T {\bf u})| \leq C \cdot max \{ d(x,y), d(y,z), d(z,x) \}$$
\end{lmm}

Because $\nabla$ and $\sigma$ are lifted from $M$, it suffices to prove each lemma on $M$.  Either connection on $M$ gives rise to horizontal spaces in the tangent bundle $TM$ and the frame bundle $\mathcal{F}M$.  Then one may endow these bundles with Riemannian metrics such that the projections to $M$ are Riemannian submersions.  We fix one such metric on $TM$ and one on $\mathcal{F}M$, and we denote the resulting norm and distance by $| \cdot |$ and $d$, as on $M$.

We first define some objects that will be referred to in the proofs of each of the three lemmas.  Consider $exp^{\sigma}$ as a map from a neighborhood of the zero section in $TM$ to $M \times M$:
$$exp^{\sigma} : (x,p) \rightarrow (x, exp^{\sigma}_x(p))$$

Then the image $exp^{\sigma}(B^{\epsilon}M)$ consists of all pairs $(x,y) \in M \times M$ with $d(x,y) < \epsilon$.  Consider $exp^{\nabla}$ the same way, and let
 \begin{eqnarray*}
V & = & (exp^{\nabla})^{-1}(exp^{\sigma}(B^{2 \epsilon}M)) \subset TM \\
U & = & (exp^{\nabla})^{-1}(exp^{\sigma}(B^{\epsilon}M)) \subset TM
\end{eqnarray*}

The map 
\begin{eqnarray*} 
f & : & V \rightarrow B^{2 \epsilon}M  \\
f & : &  (x,p) \mapsto (x, (exp^{\sigma}_x)^{-1} \circ exp^{\nabla}_x(p))
\end{eqnarray*}

is a diffeomorphism because $B(x, 2 \epsilon)$ is a normal neighborhood for all $x$.  Because $U \subset \overline{U} \subset V$, and $\overline{U}$ is compact, there exist $C_1 , C_2 > 0$ such that
\begin{eqnarray*}
d(p, q) & \leq & C_1 \cdot d(f(p),f(q)) \qquad  \mbox{for all}\ p,q \in U \\
|f_*(u)| & \leq & C_2 |u|  \qquad \mbox{for all}\ u \in TU 
\end{eqnarray*}

For $x \in M$, let 
$$f_x = (exp^{\sigma}_x)^{-1} \circ exp^{\nabla}_x$$

the restriction of $f$ to the fiber over $x$.

\begin{Pf}(of \ref{one.nabla.vs.sigma})

The map 
$$ exp^{\sigma} : B^{2 \epsilon}M \rightarrow M \times M$$

is a diffeomorphism onto its image.  There exists $C_3 >0$ such that 
$$ | (exp^{\sigma}_x)_* (u)| \leq C_3 |u| $$

for all $x \in M$ and $u$ tangent to $B({\bf 0},\epsilon) \subset T_xM$.

Let $C = C_1 C_2 C_3$.  Suppose $x,y \in M$ with $d(x,y) < \epsilon$.  We now work in $U \cap T_xM$.  Let $p = {\bf 0}$ and $q = (exp^{\nabla}_x)^{-1}(y)$.  Let $\alpha : [0,1] \rightarrow T_xM$ be the line from $p$ to $q$.  Then 
$$|\alpha^{\prime}(0)| = d(p,q) \leq C_1 \cdot d({\bf 0}, f_x(q)) = C_1 \cdot d(x,y) $$

Let $\beta$ be the image curve $f_x \circ \alpha$.  Then for any $t \in [0,1]$,
$$ | \beta^{\prime}(t)| \leq C_2 \cdot |\alpha^{\prime}(t)| = C_2 \cdot |\alpha^{\prime}(0)| \leq C_2 C_1 d(x,y)$$

The $\nabla$-geodesic $\gamma$ from $x$ to $y$ is the image $exp^{\sigma}_x \circ \beta$.  The lemma is now proved with
$$ l(\gamma) \leq sup_t |(exp^{\sigma}_x)_*(\beta^{\prime}(t))| \leq C_3 \cdot sup_t |\beta^{\prime}(t)| \leq C \cdot d(x,y)$$
\end{Pf}

\begin{Pf} (of \ref{two.nabla.vs.sigma})
\nopagebreak[3]

Define
$$ \widehat{V} = \{ (x,{\bf p},t) : (x,{\bf p}) \in T^1M, t \in \R, \ \mbox{and}\ (x,t{\bf p}) \in V \}$$

Define $\widehat{\overline{U}}$ and $\widehat{U}$ similarly.  The function
$$ Q(x,{\bf p},t) = \frac{| f_x(t {\bf p}) - t {\bf p}|}{t^2}$$

is obviously continuous on $\widehat{V} \backslash (T^1M \times \{ 0 \})$.  For each $x \in M$, the map $f_x$ is twice differentiable at ${\bf 0}$; moreover, $f_x({\bf 0}) = {\bf 0}$ and $(f_x)_{*\bf 0} = Id$.  Therefore, for any $(x,{\bf p}) \in T^1 M$, the following limit exists:
$$ \lim_{t \rightarrow 0} \frac{|f_x(t {\bf p}) - t {\bf p}|}{t^2} = \lim_{t \rightarrow 0} Q(x, {\bf p},t) = \left| D^2f_x({\bf p}, {\bf p}) \right|$$

Since $f$ is smooth, $Q$ extends to a continuous function on $\widehat{V}$.  Let $C$ be the maximum value of $Q$ on $\widehat{\overline{U}}$.  Then, for every $(x,{\bf p}) \in U$, 
$$ |f({\bf p}) -{\bf  p}| \leq C \cdot |{\bf p}|^2$$

Let ${\bf u}, {\bf v}$ be vectors in $\R^k$ equipped with the standard Euclidean inner product. Then  
$$ \angle({\bf u}, {\bf v}) \leq  \pi \cdot \left| \frac{{\bf u}}{| {\bf u}|} - \frac{{\bf v}}{|{\bf v}|} \right|$$

Indeed, consider two points $u, v$ on the unit sphere.  Let $\alpha$ be the distance between them on the unit sphere and $d$ their distance in $\R^k$.  When $0 \leq \alpha \leq \pi/2$, then $\sin \alpha \leq d$, so $\alpha/d \leq \alpha / \sin \alpha$.  When $\pi/2 \leq \alpha \leq \pi$, then $ 1 - \cos \alpha \leq d$, so $ \alpha /d \leq \alpha/ (1 - \cos \alpha)$. It is not hard to show that 
$$ \alpha / \sin \alpha \leq \pi \qquad 0 \leq \alpha \leq \pi/2$$

and 
$$ \alpha / (1 - \cos \alpha) \leq \pi \qquad \pi/2 \leq \alpha \leq \pi$$

Because each $T_xM$ equipped with the inner-product $\sigma_x$ is isometric to $\R^k$, the angle between any ${\bf p}, {\bf q} \in T_xM$ satisfies
$$ | \angle({\bf p}, {\bf q})| \leq \pi \cdot \left|\frac{{\bf p}}{|{\bf p}|} - \frac{{\bf q}}{|{\bf q}|} \right|$$

Let $x,y \in M$ with $d(x,y) < 2 \epsilon$.  Let $f_x$ be as above, and let ${\bf q} = (exp^{\nabla}_x)^{-1}(y)$.  Then we have
$$ | \angle(\gamma_{\nabla}^{\prime}(0),\gamma_{\sigma}^{\prime}(0)) | \leq \pi \cdot \left| \frac{{\bf q}}{|{\bf q}|} - \frac{f_x({\bf q})}{|f_x({\bf q})|} \right|$$

Now let $\delta \leq \epsilon$ be such that $exp^{\nabla}(B^d M \cap U)$ contains a $\delta$-neighborhood of the diagonal of $M \times M$, where $d = \theta / 2 \pi C$.  

Given $x,y \in M$, with $d(x,y) < \delta$, let ${\bf q} = (exp^{\nabla}_x)^{-1}(y)$.  Note that $(x, {\bf q}) \in U$, and $| {\bf q}| < \theta / 2 \pi C$.  Then
\begin{eqnarray*}
| \angle(\gamma_{\nabla}^{\prime}(0),\gamma_{\sigma}^{\prime}(0)) | & \leq & \pi \cdot \left| \frac{{\bf q}}{|{\bf q}|} - \frac{f_x({\bf q})}{|f_x({\bf q})|} \right|  \\
  & = &  \pi \cdot \frac{\left| {\bf q} - \frac{|{\bf q}|}{|f_x({\bf q})|} \cdot f_x({\bf q}) \right|}{|{\bf q}|}  \\
  & \leq & \pi \cdot \frac{|{\bf q} - f_x({\bf q})|}{| {\bf q}|} + \pi \cdot \frac{ \left| f_x({\bf q}) - \frac{|{\bf q}|}{|f_x({\bf q})|} \cdot f_x({\bf q}) \right|}{|{\bf q}|}  \\
  & \leq & 2  \pi \cdot \frac{|f_x({\bf q}) - {\bf q}|}{|{\bf q}|}  \\
  & \leq  & 2 \pi C \cdot |{\bf q}| \\ 
  & <  &  \theta 
\end{eqnarray*}
\end{Pf}

\begin{Pf} (of \ref{three.nabla.vs.sigma})

For any $x,y \in M$ with $d(x,y) < \epsilon$, let $t \in \R$ and ${\bf q} \in T^1_xM$ be such that 
$$ exp^{\nabla}_x(t {\bf q}) = y$$

Note that $(x, t {\bf q}) \in U$.  Let $C_1$ be as above, so
$$ t = d({\bf 0}, t {\bf q}) \leq C_1 \cdot d({\bf 0}, f(t {\bf q})) = C_1 \cdot d(x,y)$$

Next we will define a bundle with a flow that expresses parallel transport along geodesics.  Let $\pi : \mathcal{F}M \rightarrow M$ be the natural projection.  Parallel transport acts on frames $\tau \in \mathcal{F} M$ by $(P_{\gamma}^t \tau )({\bf u}) = P_{\gamma}^t (\tau({\bf u}))$.  For any subset $K \subset \mathcal{F}M$, let 
$$ \mathcal{I}K = \{ (\tau, {\bf v}) : \tau \in K, {\bf v} \in T_{\pi(\tau)}M \cap T^1M \}$$

The fibers of $\mathcal{I} K \rightarrow K$ are compact, so $\mathcal{I} K$ is compact if $K$ is.

For ${\bf v} \in T^1M$, denote by $\gamma_{\bf v}$ the geodesic with initial vector ${\bf v}$.  Now define
\begin{eqnarray*}
\mathcal{P} & : & \R^+ \times \mathcal{I} \mathcal{F} M \rightarrow \mathcal{F}M \\
\mathcal{P} & : & (t,(\tau, {\bf v})) \mapsto P_{\gamma_{\bf v}}^t \tau
\end{eqnarray*}

This map is smooth.  Let $K_0 = \mathcal{O}M$.  There exists $D_1 > 0$ such that  
$$ d(\mathcal{P}(0,(\tau, {\bf v})), \mathcal{P}(t, (\tau, {\bf v})) = d(\tau, P_{\gamma_{\bf v}}^t \tau) \leq D_1 \cdot t$$ 

for all $(\tau, {\bf v}) \in \mathcal{I} K_0$ and $t \leq C_1 \epsilon$.  Let $K_1 = \mathcal{P}([0,C_1 \epsilon] \times \mathcal{I}K_0)$.  It is compact.  There exists $D_2 > 0$ such that
$$ d(\mathcal{P}(0,(\tau, {\bf v})), \mathcal{P}(t, (\tau, {\bf v})) = d(\tau, P_{\gamma_{\bf v}}^t \tau) \leq D_2 \cdot t$$

for all $(\tau, {\bf v}) \in \mathcal{I} K_1$ and $t \leq C_1 \epsilon$.  Let $K_2 = \mathcal{P}([0,C_1 \epsilon] \times \mathcal{I}K_1)$.  There exists $D_3 > 0$ such that
$$ d(\mathcal{P}(0,(\tau, {\bf v})), \mathcal{P}(t, (\tau, {\bf v})) = d(\tau, P_{\gamma_{\bf v}}^t \tau) \leq D_3 \cdot t$$

for all $(\tau, {\bf v}) \in \mathcal{I} K_2$ and $t \leq C_1 \epsilon$.  Let $K_3 = \mathcal{P}([0, C_1 \epsilon] \times \mathcal{I}K_2)$.  

Next we identify fibers of $K_3$ with subsets of $GL(k)$ and bound the distortion of this identification.  Let $\pi^*K_3 \rightarrow K_3$ be the pullback bundle for $\pi : K_3 \rightarrow M$; elements are $(\tau_0, \tau_1)$, where $\tau_0, \tau_1 \in K_3$ and $\pi(\tau_0) = \pi(\tau_1)$.  Define
\begin{eqnarray*}
\theta & : & \pi^*K_3 \rightarrow GL(k) \\
\theta & : & (\tau_0, \tau_1) \mapsto g \ \mbox{where} \ \tau_1 = \tau_0 \circ g^{-1}
\end{eqnarray*}

Fix any metric on $GL(k)$ giving rise to a distance $d$. There exists $A_1 > 0$ such that for all $\tau_0 \in K_3$, 
$$ d(\theta(\tau_0, \tau_1), \theta(\tau_0, \tau_2)) \leq A_1 \cdot d(\tau_1, \tau_2)$$

for all $\tau_1, \tau_2 \in \pi^{-1} ( \pi(\tau_0)) \cap K_3$.  The image of $\theta$ is compact, so there exists $A_2 > 0$ such that for all $g \in im(\theta)$ and ${\bf x } \in \R^k$, 
$$\angle(g{\bf x}, {\bf x}) \leq A_2 \cdot d(e,g)$$

Let $C = C_1 A_1 A_2 (D_1 + D_2 + D_3)$.  Let $T = \Delta(xyz)$ be a $\nabla$-geodesic triangle and $-T$ the same triangle with the reverse parametrization.  Note that $P_T^{-1} = P_{-T}$.  Let the sides of $-T$ be 
\begin{eqnarray*}
\gamma_{\bf u}(0) = x  & \qquad & \gamma_{\bf u}(t_1) = z \\
\gamma_{\bf v}(0) = z  & \qquad & \gamma_{\bf v}(t_2) = y \\
\gamma_{\bf w}(0) = y  & \qquad & \gamma_{\bf w}(t_3) = x
\end{eqnarray*}

where $|{\bf u}| = |{\bf v}| = |{\bf w}| = 1$.  Let $d = max \{ d(x,y), d(y,z), d(z,x) \}$.  The parameters
$$ t_1, t_2, t_3 \leq C_1 \cdot d < C_1 \epsilon $$

Pick any $\tau_0 \in \pi^{-1}(x) \cap K_0$.  Let
\begin{eqnarray*}
\tau_1 & = & P_{\gamma_{\bf u}}^{t_1} \tau_0 \\
\tau_2 & = & P_{\gamma_{\bf v}}^{t_2} \tau_1 \\
\tau_3 & = & P_{\gamma_{\bf w}}^{t_3} \tau_2
\end{eqnarray*}

The frames $\tau_i$ belong to $K_i$ for $i = 0,1,2,3$.  Note that $\tau_3 = P_{-T} \tau_0 = P_T^{-1} \tau_0$.

We have
\begin{eqnarray*}
d(\tau_0, P_T^{-1} \tau_0 ) & = & d(\tau_0, \tau_3) \\ 
  & \leq & d(\tau_0, \tau_1) + d( \tau_1, \tau_2) + d(\tau_2, \tau_3)   \\
  & \leq & D_1 t_1 + D_2 t_2 + D_3 t_3  \\
  & \leq & C_1 (D_1 + D_2 + D_3) d
\end{eqnarray*}

Finally, for ${\bf u} \in T_xM$
\begin{eqnarray*}
\angle_{\sigma}({\bf u}, P_T({\bf u})) & = & \angle(\tau_0^{-1}({\bf u}), \tau_0^{-1}(P_T {\bf u})) \\
& = & \angle(\tau_0^{-1}({\bf u}), \tau_3^{-1}({\bf u})) \\
& = & \angle(\tau_0^{-1}({\bf u}), \theta(\tau_0,\tau_3) \cdot \tau_0^{-1} ({\bf u})) \\
& \leq & A_2 \cdot d(e, \theta(\tau_0, \tau_3)) \\
& \leq & A_1 A_2 \cdot d(\tau_0, \tau_3) \\
& \leq & C_1 A_1 A_2 (D_1 + D_2 + D_3) d
\end{eqnarray*}
\end{Pf}

For any triangle in $\R^k$ with sides of length $a, b$, and $c$, let $\theta$ be the angle opposite the side of length $c$.  Assume $a \geq b$.  Then
\begin{eqnarray*}
c^2 & = & a^2 + b^2 - 2ab \cos \theta \\
& \geq & a^2 + b^2 - 2a^2 \cos \theta \\
& \geq & a^2 + (c-a)^2 - 2a^2 \cos \theta
\end{eqnarray*}

from which we obtain
$$ a \leq c /(1 - \cos \theta)$$

Then the proof of this lemma is straightforward using the exponential map for $\sigma$.
\begin{lmm}
\label{triangle.lengths}
Let $\epsilon > 0$ be as above, so that every $B(x, 2 \epsilon)$ is a normal neighborhood of each of its points.
There exists $C > 0$ such that, in any $\sigma$-geodesic triangle $\alpha \beta \gamma = \Delta(xyz)$ in $X$ with 
$$ \mbox{max} \{ d(x,y), d(y,z), d(z,x) \}  < \epsilon$$

with $\theta$ the angle opposite $\gamma$.
$$ max \{ l(\alpha), l(\beta) \} \leq C \cdot \frac{l(\gamma)}{1 - \cos \theta}$$

\end{lmm}

We now proceed to the proof of the proposition, following the idea of \cite{ZeLip} Lemma 14.

\begin{Pf} 
  As above, we will fix a Riemannian metric on $TM$ such that $ \pi : TM \rightarrow M$ is a Riemannian submersion with respect to $\sigma$.  We can lift this metric to a $\Gamma$-invariant metric on $TX$.  The resulting distance will in all cases be denoted by $d$.

Let $r > 0$ be given.  There exists $\eta > 0$ such that $B(x, \eta) \subset exp_x (B({\bf 0}, r))$ for all $x \in X$.  Let $0 < \epsilon < \eta/3$.  Then $\epsilon$ has the property
\begin{eqnarray}
\forall \ p \in X,\ \forall \ x \in B(p, 3 \epsilon / 2):  exp_x^{-1}(D(p, 3 \epsilon / 2)) \subset B({\bf 0}, r) \subset T_xX    \label{epsilon.prop.r}
\end{eqnarray}

We may assume that $\epsilon > 0$ has the following further properties:
\begin{eqnarray}
\forall \ p \in X & : & B(p, 2 \epsilon) \ \mbox{is a $\sigma$- and $\nabla$-normal neighborhood of each point}  \label{epsilon.prop.normal}\\
\forall \ p \in X & : & B(p, 2 \epsilon) \ \mbox{is hypersurface foliated}  \label{epsilon.prop.hyp.fol}
\end{eqnarray}
where $\epsilon$ with property (\ref{epsilon.prop.hyp.fol}) is given by Lemma \ref{hypersurface.foliated}.

Last, for any $x,y \in X$ with $d(x,y) < \epsilon$, 
\begin{eqnarray}
| \angle(\gamma^{\prime}(0), \gamma_{\sigma}^{\prime}(0)) | \leq \pi / 8   \label{epsilon.prop.angles}
\end{eqnarray}

where $\gamma$ and $\gamma_{\sigma}$ are the $\nabla$- and $\sigma$-geodesics, respectively, from $x$ to $y$. Such an $\epsilon > 0$ exists by Lemma \ref{two.nabla.vs.sigma}.

By property (\ref{epsilon.prop.normal}), the exponential map is a diffeomorphism from a neighborhood of any $(p,{\bf 0})$ onto $B(p,2 \epsilon) \times B(p, 2 \epsilon)$, so there exists 
$$ \psi_p : B(p, 3\epsilon/2) \times B(p, 2 \epsilon) \rightarrow TX$$ 

such that $exp \circ \psi_p$ is the identity on $B(p, 3 \epsilon/2) \times B(p,2 \epsilon)$.  Since $\psi_p$ is smooth, there exists $C_p$ such that 
\begin{eqnarray}
d( \psi_p(x,w), \psi_p(y,w)) \leq C_p \cdot d(x,y)  \label{psi.lipschitz}
\end{eqnarray}

for all $x,y \in D(p, \epsilon)$ and $w \in D(p, 3 \epsilon/2)$.  

The image $\psi_p( D(p,\epsilon) \times \partial D(p, 3 \epsilon / 2))$ is a compact subset of $T D(p, \epsilon ) \backslash Z$, where $Z$ is the zero section of $T D(p,\epsilon)$.  Here we have the smooth map
\begin{eqnarray*}
\varphi_p & : & T D(p, \epsilon) \backslash Z \rightarrow T^1 D(p, \epsilon) \\
\varphi_p & : & (p, {\bf u}) \mapsto (p, {\bf u} / |{\bf u}|)
\end{eqnarray*}

There exists $D_p > 0$ such that 
\begin{eqnarray}
d(\varphi_p (x, {\bf u}), \varphi_p(y,{\bf v}) \leq D_p \cdot d((x, {\bf u}), (y, {\bf v} ))  \label{varphi.lipschitz}
\end{eqnarray}

for all $(x,{\bf u}), (y,{\bf v}) \in \psi_p( D(p, \epsilon) \times \partial D(p, 3 \epsilon / 2))$.  

Let $\overline{p}_1, \ldots, \overline{p}_N$ be such that $B(\overline{p}_1, \epsilon), \ldots, B(\overline{p}_N, \epsilon)$ is a finite cover of $M$.  For any lifts $p_1, \ldots, p_N$ of these points to $X$, these $\epsilon$-balls lift to balls $B(p_i, \epsilon)$ whose translates $B(\gamma p_i, \epsilon)$ cover $X$.  Let $C_i$ be the Lipschitz constant for $\psi_{p_i}$ on $D(p_i, \epsilon) \times D(p_i, 3 \epsilon / 2)$, for $i = 1, \ldots, N$, and let $C_0$ be the maximum of these.  Let $D_i$ be the Lipschitz constant for $\varphi_{p_i}$ on $\psi_{p_i} ( D(p_i, \epsilon) \times D(p_i, 3 \epsilon / 2))$, for $i = 1, \ldots, N$, and let $D_0$ be the maximum of these.  For any $p = \gamma p_i$, the map $\psi_p$ has the Lipschitz property (\ref{psi.lipschitz}) with Lipschitz constant $C_0$ on $D(p,\epsilon) \times D(p,3 \epsilon /2)$, and $\varphi_p$ has property (\ref{varphi.lipschitz}) with Lipschitz constant $D_0$ on $\psi_p( D(p, \epsilon) \times D(p,3 \epsilon / 2))$.  We assume $\delta_0 < \epsilon/2$.  Let $\delta_0$ be the Lebesgue number of the covering by the balls $B( \gamma p_i, \epsilon)$ of $X$.   

Let $\eta > 0$ be such that $exp$ is defined on $B^{\eta} X$ and $exp(B^{\eta}X)$ is contained in the $\delta_0$-neighborhood of the diagonal in $X \times X$.  By Lemma \ref{normally.foliated}, there exists $\delta_1 > 0$ such that 
\begin{eqnarray}
\forall \ x \in X & : & B(x, \delta_1) \ \mbox{is normally foliated based at} \ x   \label{delta1.prop.one} \\
\forall \ x \in X, H \subset T_xX  & : & \varphi_{H,x}^{-1} (B(x, \delta_1)) \subseteq B({\bf 0}, \eta)  \subseteq H \times \R \label{delta1.prop.two}
\end{eqnarray}

Finally, take $\delta_2$ such that, 
\begin{eqnarray}
\forall y,z \in B(x, \delta_2),\  \mbox{the $\sigma$-geodesic from $y$ to $z$ is contained in} \ B(x,2\epsilon) \label{delta2.prop}  
\end{eqnarray}

Now let $\delta = min\{ \delta_0, \delta_1, \delta_2 \}$.

Let $C_1, C_2 > 0$ be the constants given by Lemmas \ref{three.nabla.vs.sigma} and \ref{triangle.lengths}, respectively, with the $\epsilon$ we have chosen above.  By property (\ref{epsilon.prop.normal}), this $\epsilon$ satisfies the hypotheses of the lemmas.

Let $(X^{\prime}, f)$ be any radius-$r$ codimension-one geodesic lamination on $X$.  Let $x,y \in X^{\prime}$ with $d(x,y) < \delta$.  We may assume $\mathcal{L}_x \cap \mathcal{L}_y = \emptyset$.  This intersection is open in each plaque (see Definition \ref{defn.geod.lamn}); if it were nonempty, then the geodesic $\gamma$ from $x$ to $y$ would be contained in $\mathcal{L}_x \cap \mathcal{L}_y$, and $P_{\gamma} f(x) = f(y)$.  Let $p = \gamma p_i \in X$ be such that $B(x, \delta) \subset B(p, \epsilon)$.  

Now $\psi_p$ satisfies 
$$ d( \psi_p(z,w), \psi_p(y,w)) \leq C_0 \cdot d(z,y)$$

for all $z \in D(p, \epsilon)$ and $w \in D(p, 3 \epsilon / 2)$.  By properties (\ref{epsilon.prop.r}) and (\ref{epsilon.prop.normal}) of $\epsilon$, the intersection $\mathcal{L}_x \cap D(p, 3 \epsilon/2)$ is contained in the convex hull of $S_x = \mathcal{L}_x \cap \partial D(p, 3 \epsilon / 2)$.  For any $z \in B(p, \epsilon)$, let 
$$R(z) = \varphi_p \circ \psi_p ( \{ z \} \times S_x)$$

The set $R(x)= f(x) \cap T^1_x X$.  The set $R(y)$ consists of all ${\bf u}$ in $T^1_y X$ such that $exp_y (t {\bf u}) \in S_x$ for some positive $t$.  The assumption that $y \notin \mathcal{L}_x$ implies $ - R(y) \cap R(y) = \emptyset$.

By property (\ref{delta1.prop.one}) of $\delta$, the ball $B(x, \delta)$ is normally foliated based at $x$.  Choose an orientation of $f(x)$.  There exist ${\bf q} \in f(x)$ and $t > 0$ such that 
$$ y = exp_{\gamma_{{\bf q}(1)}} ( t N (P^1_{\gamma_{\bf q}} f(x)))$$

Let $z = \gamma_{\bf q}(1) = exp_x({\bf q})$.  By property (\ref{delta1.prop.two}) of $\delta$, the vector ${\bf q} \in B({\bf 0}, \eta)$, so $z \in \mathcal{L}_x \cap B(x, \delta_0) \subset \mathcal{L}_x \cap B(p, \epsilon)$.  Let $\beta$ be the geodesic from $z$ to $x$.  Denote by $f(z)$ the hypersurface $P_{- \beta} f(x) = T_z \mathcal{L}_x$.  Let $\alpha$ be the geodesic from $z$ to $y$ with $\alpha^{\prime}(0) = N(f(z))$.  

Let $\gamma$ be the geodesic from $x$ to $y$, and let $T$ be the geodesic triangle $\alpha \beta \gamma$, where $\alpha$ is traversed from $y$ to $z$.  Then 
$$P_{\gamma} f(x) = P_T P_{\alpha} f(z)$$  

Because $x,y,z \in B(x, \delta_0)$ and $\delta_0 < \epsilon / 2$, each of the distances $d(x,y), d(y,z), d(z,x) < \epsilon$.  By Lemma \ref{three.nabla.vs.sigma}, for any ${\bf u} \in P_{\alpha} f(z)$, 
$$ | \angle( {\bf u}, P_T {\bf u})| \leq C_1 \cdot max \{ d(x,y), d(y,z), d(z,x) \}$$

Let $\alpha_{\sigma}$ be the $\sigma$-geodesic from $z$ to $y$, $\beta_{\sigma}$ the $\sigma$-geodesic from $z$ to $x$, and $\gamma_{\sigma}$ the $\sigma$-geodesic from $x$ to $y$.  By property (\ref{epsilon.prop.normal}) of $\epsilon$, the length $l(\gamma_{\sigma}) = d(x,y)$.  By property (\ref{epsilon.prop.angles}) of $\epsilon$, the angles 
$$ | \angle ( \alpha^{\prime}(0), \alpha_{\sigma}^{\prime}(0))|, \ | \angle(\beta^{\prime}(0), \beta_{\sigma}^{\prime}(0))| \leq \pi /8$$

Therefore, 
$$ \pi / 4 = |\angle(\alpha^{\prime}(0), \beta^{\prime}(0))| - \pi / 4 \leq |\angle(\alpha_{\sigma}^{\prime}(0), \beta_{\sigma}^{\prime}(0)| \leq |\angle(\alpha^{\prime}(0), \beta^{\prime}(0))| + \pi / 4 = 3 \pi /4$$

Because $\mbox{max} \{ d(x,y), d(y,z), d(z,x) \} < \epsilon$, Lemma \ref{triangle.lengths} applies to the triangle $\alpha_{\sigma} \beta_\sigma \gamma_{\sigma}$, and there exists $C_2 > 0$ such that 
$$ l(\alpha_{\sigma}), l(\beta_{\sigma}) \leq \frac{2 C_2 \cdot l(\gamma_{\sigma})}{2- \sqrt{2}} = \frac{ 2 C_2 \cdot d(x,y)}{2 - \sqrt{2}} < 4 C_2 \cdot d(x,y)$$  

so $max \{ d(x,y), d(y,z), d(z,x) \} \leq 4 C_2 \cdot d(x,y)$.  Then we have, for any ${\bf u} \in P_{\alpha} f(z)$, 
$$ | \angle( {\bf u}, P_T {\bf u})| \leq 4 C_1 C_2 \cdot d(x,y) $$

Therefore, 
$$ |\angle (P_{\alpha} f(z), P_{\gamma} f(x)) | \leq 4 C_1 C_2 \cdot d(x,y)$$

By property (\ref{epsilon.prop.hyp.fol}) of $\epsilon$, the ball $B(p, 2 \epsilon)$ is hypersurface foliated based at $z$, so $exp_y (P_{\alpha} f(z))$ does not intersect $\mathcal{L}_x \cap D(x, 3 \epsilon /2)$.  In particular, $P_{\alpha} f(z) \cap T^1_y X$ is contained in the annulus with boundary $- R(y) \cup R(y)$.  The condition $\mathcal{L}_x \cap \mathcal{L}_y = \emptyset$ implies that $f(y) \cap T^1_yX$ is also contained in this annulus.  

We will show this annulus is narrow using the Lipschitz properties (\ref{psi.lipschitz}) and (\ref{varphi.lipschitz}).  Let $d_H$ be the Hausdorff distance on compact sets in $TX$.  
\begin{eqnarray*}
d_H(- R(y), R(y)) & \leq & d_H(- R(y), R(x)) + d_H (R(x), R(y)) \\
& = & d_H(- R(y), - R(x)) + d_H(R(x), R(y)) \\
& \leq & 2 D_0 \cdot d_H( \psi_p (\{ x \} \times S_x), \psi_p ( \{ y \} \times S_x)) \\
& \leq & 2 D_0 C_0 \cdot d(x,y)
\end{eqnarray*}

Then 
\begin{eqnarray*}
| \angle(P_{\alpha} f(z), f(y)) | & \leq & \pi \cdot d_H( P_{\alpha} f(z) \cap T^1_yX, f(y) \cap T^1_y X)  \\
  & \leq & 2 \pi \cdot d_H(- R(y), R(y)) \\
& \leq & 4 \pi D_0 C_0 \cdot d(x,y)
\end{eqnarray*}

Finally
\begin{eqnarray*}
| \angle (P_{\gamma} f(x), f(y)) | & \leq & | \angle(P_{\gamma} f(x), P_{\alpha} f(z)) | + | \angle(P_{\alpha} f(z), f(y))| \\
& \leq & 4 (C_1 C_2 + \pi C_0 D_0) \cdot d(x,y)
\end{eqnarray*}
\end{Pf}

\end{document}